\documentclass{amsart}
\newcommand{\note}{\noindent {\bf Notation. }}
\newcommand{\remark}{\noindent {\bf Remark. }}
\newcommand{\ws}{\hspace{4pt}}
\newtheorem{theorem}{Theorem}

\newtheorem{corollary}{Corollary}
\newtheorem{lemma}{Lemma}
\newtheorem{defi}{Definition}
\usepackage[]{amssymb,amsopn}
\usepackage[]{amsmath}
\usepackage[]{eucal}
\usepackage[]{latexsym}

\begin{document}

\title{Potential Theory and Quadratic Programming}
\author{\'A. P. Horv\'ath }


\begin{abstract}  We extend the notion of some energy-type expressions based on two sets, developed in the abstract potential theory. We also give the discretized version of the quantities defined, similar to Chebyshev constant. This extension allows to apply the potential-theoretic results to infinite quadratic programming problems. Together with a cutting plane algorithm, the Chebyshev-constant method ensures that under certain conditions, the infinite problem can be reduced to semi-infinite or to finite problems.\end{abstract}
\maketitle

\section{Introduction}

A rather general setting of abstract (linear) potential theory was developed by several authors mainly in the '60-s, e.g. Frostman \cite{f}, Choquet \cite{ch}, Fuglede \cite{fu}, Ohtsuka \cite{o}, etc. The basic problem was to find the Wiener energy of a (compact) set in a locally compact Hausdorff space that is
$$w(K)=\inf_{\mu\in \mathcal{M}_1(K)}\int_{X\times X}k(x,y)d\mu(x)\times \mu(y),$$
where $\mathcal{M}_1(K)$ is the set of probability measures on $K$ and $k$ is a lower semicontinuous kernel function on $X\times X$. Several energy-type expressions were defined and their relations were investigated by several authors. For instance it is proved (\cite{f}, see also \cite{fana}) that under some assumptions, the two expressions $w(K)$ and $q(K)=\inf_{\mu\in \mathcal{M}_1(K)}\sup_{y\in K}\int_Xk(x,y)d\mu(x)$ coincide. This immediately led to study problems with two different sets (sometimes from two different spaces), cf. \cite{fu1}, \cite{o1}. This two sets method proved to be useful tool to examine the so-called rendezvous numbers \cite{fare}. In this context the previous definition is modified as $q(K,L)=\inf_{\mu\in \mathcal{M}_1(K)}\sup_{y\in L}\int_Xk(x,y)d\mu(x)$, which can be expressed by Dirac measures concentrated at the points of $L$. Changing the set of measures at the second variable, we arrive to $p$-capacity as it is defined in the book of Adams and Hedberg \cite{ah}.  It is $C_p(K)^{-\frac{1}{p}}=\inf_{\mu\in \mathcal{M}_1(K)}\sup_{\{fd\nu : f\geq 0, \|f\|_{p,\nu}\leq 1\}}\int_{\mathbb{R}^n}\int_Xk(x,y)d\mu(x)f(y)d\nu(y)$. The problem is the following. Replacing the original sets of measures (probability measures, Dirac measures), in the two variables of the minimax problem with some other sets of measures, under what conditions can we state some coincidence with a Wiener energy-type expression when $K=L$, say. It turns out that the original construction works, because the second set of measures is the extremal set of the first one. This observation orients our attention to the problem of quadratic programming.

The relations between potential theory and linear programming was discovered also in the '60-s, and the infinite linear programming was developed by Ohtsuka, Yoshida, etc. (cf. e.g. \cite{o2}, \cite{o3}, \cite{yo}, \cite{ny}). For further results see \cite{rh}. A nice application of quadratic programming to computation of capacity is given by Rajon, Ransford and Rostand, cf. \cite{rrr} and the references therein. Our field of interest is infinite quadratic programming which is related to infinite linear programming and to potential theory as well. For drafting the problem see \cite{cw}, \cite{w}, \cite{wfl}, \cite{wwt}. The aim of our investigation is to convert infinite problems to a limit of semi-infinite or finite problems by mixed - potential-theoretic and quadratic programming - methods. Besides the cutting plane method, our main tool is the so-called Chebyshev constant, which is the discrete version of potential.

This paper is organized as follows. Section 2 describes the framework and introduces terminology used throughout this paper. Section 3 extends the definition of some energy-type expression and investigates the basic properties of the quantities defined. Section 4 introduces the quadratic programming problem connected with the energy defined in Section 3 and gives discretization methods. Finally in the Appendix there are several examples on sets of measures and kernel functions which satisfy the required properties.

\section{Definitions and Notation}

$X$, $Y$ are locally compact Hausdorff spaces, the kernel $k: X\times Y \to [0,\infty]$ is lower semicontinuous (l.s.c.).
$$\mathcal{M}(X):=\left\{\mu : \mu \ws \mbox{is a positive, regular Radon measure on $X$}, \ws \mu(X)<\infty\right\},$$
and the set of probability measures from $\mathcal{M}(X)$ is $\mathcal{M}_1(X):=\left\{\mu \in \mathcal{M}(X) : \|\mu\|=1\right\}$.
Let $H\subset X$, $L\subset Y$ arbitrary. We say that $\mu\in \mathcal{M}(X)$ is {\it concentrated} on the set $H$ if each compact set intersects the complement of $H$ in a set of zero (outer-) measure; or equivalently, if $H$ is $\mu$-measurable and $\mu=\mu|_H$ cf. \cite[p. 146.]{fu}. (We also denote by $\mu$ the extension of $\mu$ to an outer measure.) We say that $\mu$ is {\it supported} on the set $H$ if $\mathrm{supp}\mu\subset H$.
$$\mathcal{M}_1(H):=\left\{\mu \in \mathcal{M}_1(X) : \mu \ws \mbox{is concentrated on}\ws H\right\}.$$
Similarly $\mathcal{M}(H):=\left\{\mu \in \mathcal{M}(X) : \mu \ws \mbox{is concentrated on}\ws H\right\}.$
$\mathcal{R}(H)$, $\mathcal{S}(L)$ are subsets of $\mathcal{M}(H)$ and $\mathcal{M}(L)$ defined suitably at the occurrences.

The mutual energy of $\mu \in \mathcal{R}(H)$ and $\nu \in \mathcal{S}(L)$ is
$$E(\mu,\nu):= \int_{X\times Y} k(x,y)d\mu(x)\times\nu(y).$$

Denoting by $C_c(X\times Y)$ the continuous, compactly supported real valued functions on $X\times Y$, the upper integral of a positive l.s.c. function is
$$\int_{X\times Y}k(x,y)d\mu(x)\times\nu(y):=\sup_{0\leq h(x,y)\leq k(x,y)\atop h\in C_c(X\times Y)}\int_{X\times Y}h(x,y)d\mu(x)\times\nu(y).$$
We use the integral of a l.s.c. positive function in upper integral sense. If the function is continuous it coincides with the usual integral. We have the following important properties (see eg. \cite[Lemma 22.5]{d})
$$E(\mu,\nu)=\int_Y \int_X k(x,y)d\mu(x)d\nu(y)=\int_X\int_Y k(x,y)d\nu(y)d\mu(x).$$
If the function is l.s.c., positive and $\mu\times \nu$-integrable the equation above fulfils in ordinary sense as well (cf, \cite[Corollary 22.14]{d}).
Moreover if $k$ is symmetric on $X\times X$
$$E(\mu,\nu)=E(\nu,\mu).$$

\medskip

Let us introduce some further notation. Let $k: X\times Y \to [0,\infty]$ as above and $\mu\in \mathcal{M}(H)$ be arbitrary. The potential function of $\mu\in \mathcal{M}(H)$ is
$$U^{\mu}(y)=\int_X k(x,y)d\mu(x),$$
where $y\in Y$. If $k: X\times X \to [0,\infty]$, $\mathcal{R}(H)\subset \mathcal{M}(H)$, let $w(\mathcal{R}(H))$ be the energy of $H\subset X$ with respect to $\mathcal{R}(H)$
$$w(\mathcal{R}(H)):=\inf_{\mu\in\mathcal{R}(H)}E(\mu),$$
where $E(\mu):=E(\mu,\mu)$.
We say that a positive kernel satisfies the {\it Frostman's maximum principle} (cf. \cite[page 150]{fu}) if for every measure $\mu \in \mathcal{M}(X)$ of compact support
$$\sup_{x\in X}U^{\mu}(x)=\sup_{x\in \mathrm{supp}\mu}U^{\mu}(x).$$

The topology on $\mathcal{M}(H)$ is the vague topology which is the locally convex topology generated by the family of seminorms $\{\|\mu\|_f : f\in C_c(X)\}$, where $\|\mu\|_f=|\int_Xfd\mu|$. When $K$ is a compact subset of $H$ we write $K\subset\subset H$. When $K\subset\subset X$, this topology on $\mathcal{M}(K)$ coincides with the $w^*$-topology determined by $C(K)$. Let $\mathcal{K}(H)$ be the filtering family or directed set of all compact subsets of $H$ (directed upwards).

Next lemma appeared in several works, cf. e.g. \cite{d}, \cite{fu}, \cite{ho1}, \cite{fare}.

\medskip

\begin{lemma} Let $H\subset X$, $L\subset Y$. Equipping $\mathcal{M}(H)$ and $\mathcal{M}(L)$ with the vague topology the following functions $F$ are lower semicontinuous

\noindent {\rm (a)}  $F: \mathcal{M}(H)\times L \to \overline{\mathbb{R}}_+ , \ws (\mu,y) \mapsto U^{\mu}(y),$

\noindent {\rm (b)}  $F: \mathcal{M}(H)\times \mathcal{M}(L) \to \overline{\mathbb{R}}_+ , \ws (\mu,\nu) \mapsto E(\mu,\nu),$

\noindent {\rm (c)}  $F: \mathcal{M}(H) \to \overline{\mathbb{R}}_+ , \ws \mu \mapsto \sup_{y\in L} U^{\mu}(y).$
\end{lemma}

\section{Energy-type Expressions}

In this section we extend the definitions of some energy-type quantities and analyze their relations. Instead of the sets of probability measures and Dirac measures on different sets we use more general sets of measures on different sets again. As it was pointed out in \cite{fare} involving two sets allows to give these energy-type quantities (for any sets) as a supremum/infimum of similar quantities with respect to compact sets. Some applications are presented, see the examples after Lemma 4 and Lemma 5. The extension of sets of measures allows to apply the results of this section to some problems of infinite quadratic programming.

\subsection{Energy with respect to Sets of Measures}
The definitions below are the extensions of the definitions of \cite{o1} and \cite{fare}.

\begin{defi} Let $H \subset X$, $L\subset Y$ and $\mathcal{R}(H)\subset \mathcal{M}(H)$, $\mathcal{S}(L)\subset \mathcal{M}(L)$.
$$q(\mathcal{R}(H),\mathcal{S}(L)):=\inf_{\mu \in \mathcal{R}(H)}\sup_{\nu \in \mathcal{S}(L)}E(\mu,\nu)=\inf_{\mu \in \mathcal{R}(H)}\sup_{\nu \in \mathcal{S}(L)}\int_Y \int_X k(x,y)d\mu(x)d\nu(y),$$
$$\underline{q}(\mathcal{R}(H),\mathcal{S}(L)):=\sup_{\mu \in \mathcal{R}(H)}\inf_{\nu \in \mathcal{S}(L)}E(\mu,\nu)=\sup_{\mu \in \mathcal{R}(H)}\inf_{\nu \in \mathcal{S}(L)}\int_Y \int_X k(x,y)d\mu(x)d\nu(y).$$\end{defi}

\medskip

\begin{theorem}
Let $k$ be a positive, symmetric l.s.c. kernel on $X\times X$. Let $H, L \subset X$ two subsets of $X$, $\mathcal{R}(H)$, $\mathcal{S}(H)$ are convex subsets of $\mathcal{M}(H)$, $\tilde{\mathcal{R}}(L)$, $\tilde{\mathcal{S}}(L)$ are $w^{*}$-compact convex subsets of $\mathcal{M}(L)$. Moreover let us suppose that
\begin{equation}  \inf_{\nu\in\tilde{\mathcal{R}}(L)}\int_L U^{\mu}d\nu=\inf_{\nu\in\tilde{\mathcal{S}}(L)}\int_L U^{\mu}d\nu   \ws \ws \forall \mu \in \mathcal{M}(H),\end{equation}
\begin{equation}  \sup_{\nu\in\mathcal{R}(H)}\int_H U^{\mu}d\nu=\sup_{\nu\in\mathcal{S}(H)}\int_H U^{\mu}d\nu   \ws \ws \forall \mu \in \mathcal{M}(L).\end{equation}
Then
\begin{equation}\underline{q}\left(\mathcal{R}(H),\tilde{\mathcal{S}}(L)\right)=q\left(\tilde{\mathcal{R}}(L),\mathcal{S}(H)\right).\end{equation}\end{theorem}

To prove the theorem we need the following lemmas.

\medskip

\lemma
 As above let $k$ is a l.s.c. symmetric kernel on $X\times X$, $H, L \subset X$, $\mathcal{R}(H), \mathcal{S}(H)\subset \mathcal{M}(H)$, $\tilde{\mathcal{R}}(L),\tilde{\mathcal{S}}(L)\subset \mathcal{M}(L)$ such that
\begin{equation}  \inf_{\tilde{\mathcal{S}}(L)}\int_L U^{\mu}d\nu\leq\inf_{\tilde{\mathcal{R}}(L)}\int_L U^{\mu}d\nu   \ws \ws \forall \mu \in \mathcal{M}(H),\end{equation}
\begin{equation}  \sup_{\mathcal{R}(H)}\int_H U^{\mu}d\nu\leq\sup_{\mathcal{S}(H)}\int_H U^{\mu}d\nu   \ws \ws \forall \mu \in \mathcal{M}(L).\end{equation}
Then
\begin{equation}\underline{q}\left(\mathcal{R}(H),\tilde{\mathcal{S}}(L)\right)\leq q\left(\tilde{\mathcal{R}}(L),\mathcal{S}(H)\right).\end{equation}

\proof
By (4) for all $\mu \in \mathcal{R}(H)$ and $\sigma \in \tilde{\mathcal{R}}(L)$
$$\inf_{\nu \in \tilde{\mathcal{S}}(L)}\int_L \int_H  k(x,y)d\mu(x) d\nu(y)\leq \int_L U^{\mu}d\sigma=\int_H \int_L k(x,y)d\sigma(y)d\mu(x)$$
$$\leq \sup_{\lambda \in \mathcal{S}(H)}\int_H \int_L k(x,y)d\sigma(y)d\lambda(x)=\sup_{\lambda \in \mathcal{S}(H)}\int_H \int_L k(x,y)d\sigma(x)d\lambda(y).$$
Here we changed the order of variables applied (5) and the symmetry of the kernel. Taking infimum and supremum over $\sigma \in \tilde{\mathcal{R}}(L)$ and  $\mu \in \mathcal{R}(H)$ respectively we obtain the inequality of the lemma.

\medskip

\lemma\cite[Theorem 5.2]{fare} or \cite[Theorem 2.4.1]{ah}

Let $A$ be a compact convex subset of the Hausdorff topological vector space $U$ and $B$ be a convex subset of the linear space $V$. Let $f: A\times B\to (-\infty,\infty]$ be l.s.c. on $A$ for fixed $y\in B$, and assume that $f$ is convex in the first and concave in the second variable. Then
$$\sup_{y\in B}\inf_{x\in A}f(x,y)=\inf_{x\in A}\sup_{y\in B}f(x,y).$$

\medskip

\proof (of Theorem 1)
By Lemma 2
$$\underline{q}\left(\mathcal{R}(H),\tilde{\mathcal{S}}(L)\right)\leq q\left(\tilde{\mathcal{R}}(L),\mathcal{S}(H)\right).$$
Since the $w^{*}$-topology is Hausdorff, $\tilde{\mathcal{R}}(L)$ is $w^{*}$-compact convex, and $\mathcal{S}(H)$ is convex furthermore $f(\sigma,\lambda)=E(\sigma,\lambda)$ is linear in both of the variables, we can apply Lemma 3, (and then Lemma 2 and Lemma 3 again) that is

$$q\left(\tilde{\mathcal{R}}(L),\mathcal{S}(H)\right)=\inf_{\sigma \in \tilde{\mathcal{R}}(L)}\sup_{\lambda \in \mathcal{S}(H)}\int_H \int_L k(x,y)d\sigma(x)d\lambda(y)$$ $$=\sup_{\lambda \in \mathcal{S}(H)}\inf_{\sigma \in \tilde{\mathcal{R}}(L)}\int_H \int_L k(x,y)d\sigma(x)d\lambda(y)=\sup_{\lambda \in \mathcal{S}(H)}\inf_{\sigma \in \tilde{\mathcal{R}}(L)}\int_L \int_H k(x,y)d\lambda(x)d\sigma(y)$$ $$=\underline{q}\left(\mathcal{S}(H),\tilde{\mathcal{R}}(L)\right)\leq q\left(\tilde{\mathcal{S}}(L),\mathcal{R}(H)\right)=\inf_{\nu\in \tilde{\mathcal{S}}(L)}\sup_{\mu\in \mathcal{R}(H)}\int_H \int_L k(x,y)d\nu(x)d\mu(y)$$ $$=
\sup_{\mu \in \mathcal{R}(H)}\inf_{\nu\in \tilde{\mathcal{S}}(L)}\int_H \int_L k(x,y)d\nu(x)d\mu(y)=\sup_{\mu \in \mathcal{R}(H)}\inf_{\nu\in \tilde{\mathcal{S}}(L)}\int_L \int_H k(x,y)d\mu(x)d\nu(y)$$ $$=
\underline{q}\left(\mathcal{R}(H),\tilde{\mathcal{S}}(L)\right).$$

\medskip

\remark With $\mathcal{R}(H)=\mathcal{S}(H)$ convex and $\tilde{\mathcal{R}}(L)=\tilde{\mathcal{S}}(L)$ compact, convex, Theorem 1 coincides with Lemma 3.

\medskip

\note
Let $\mathcal{R}(H)$ is $w^{*}$-compact convex and $\mathcal{S}(H)=\mathrm{Ex}\mathcal{R}(H)$, the extreme points of  $\mathcal{R}(H)$, that is $\sigma \in \mathrm{Ex}\mathcal{R}(H)$ if for any $\mu, \nu \in \mathcal{R}(H)$, $\alpha \in (0,1)$, $\sigma =\alpha\mu + (1-\alpha)\nu$, it necessarily follows that $\sigma=\mu=\nu$ .
We introduce the notation
$$\underline{q}\left(\mathcal{R}(H)\right):=\underline{q}\left(\mathcal{R}(H),\mathrm{Ex}\mathcal{R}(H)\right), \ws \mbox{and} \ws q\left(\mathcal{R}(H)\right):= q\left(\mathcal{R}(H),\mathrm{Ex}\mathcal{R}(H)\right).$$

\medskip

\begin{corollary}

With the notation of Theorem 1, let $H\subset\subset X$, $L \subset\subset Y$, both $\mathcal{R}(H)$ and $\tilde{\mathcal{R}}(L)$ be $w^{*}$-compact convex sets of measures, and let  $\mathcal{S}(H)$ and $\tilde{\mathcal{S}}(L)$ be the extreme points of  $\mathcal{R}(H)$ and $\tilde{\mathcal{R}}(L)$ respectively. Then
$$\underline{q}\left(\mathcal{R}(H),\tilde{\mathcal{S}}(L)\right)=q\left(\tilde{\mathcal{R}}(L),\mathcal{S}(H)\right).$$
In particular if $H=L$, $\mathcal{R}(H)=\tilde{\mathcal{R}}(H)$, then
$$\underline{q}\left(\mathcal{R}(H)\right)=q\left(\mathcal{R}(H)\right).$$\end{corollary}

\proof
First we need that (1) and (2) are satisfied. We prove (2), (1) is similar. Let $\mathrm{co}\mathcal{S}(H)$ be the set of the convex combination of the elements of $\mathcal{S}(H)$. Since $\mathcal{S}(H)\subset \mathrm{co}\mathcal{S}(H)$,
$$B:=\sup_{\nu\in\mathcal{S}(H)}\int_H U^{\mu}d\nu \leq \sup_{\nu \in \mathrm{co}\mathcal{S}(H)}\int_H U^{\mu}d\nu=\sup_{\{\alpha_i\}, \alpha_i\geq 0, \sum \alpha_i =1\atop \nu_i\in \mathcal{S}(H)}\int_H U^{\mu}d\sum \alpha_i\nu_i$$ $$ \leq \sup_{\{\alpha_i\}, \alpha_i\geq 0, \sum \alpha_i =1}\sum \alpha_i\sup_{\nu_i\in \mathcal{S}(H)}\int_H U^{\mu}d\nu_i=B.$$
That is the supremum on the elements of $\mathcal{S}(H)$ coincides with the supremum on the elements of $\mathrm{co}\mathcal{S}(H)$. According to the Krein-Milman theorem (see \cite[Chapter II, 10.4]{sch}), $\mathrm{co}\mathcal{S}(H)$ is a $w^{*}$-dense subset of $\mathcal{R}(H)$, and so the supremum on $\mathcal{R}(H)$ coincides with the supremum on $\mathrm{co}\mathcal{S}(H)$. Indeed let $A:=\sup_{\nu\in\mathcal{R}(H)}\int_H U^{\mu}d\nu$ be finite and $\varepsilon>0$ be arbitrary. (The $A=\infty$-case is similar.)  Let $\lambda \in \mathcal{R}(H)$ such that $\int_H U^{\mu}d\lambda>A-\varepsilon$. Since $U^{\mu}$ is l.s.c. by the definition of upper integral there is a continuous function $f$ with compact support such that $ U^{\mu}\geq f\geq 0$ and $A\geq \int_H U^{\mu}d\lambda\geq\int_H f d\lambda>\int_H U^{\mu}d\lambda-\varepsilon>A-2\varepsilon$. Since $\mathrm{co}\mathcal{S}(H)$ is $w^{*}$-dense in $\mathcal{R}(H)$, there is a $\sigma \in \mathrm{co}\mathcal{S}(H)$ such that $A+\varepsilon>\int_H f d\sigma>A-3\varepsilon$. That is with $\sigma \in \mathrm{co}\mathcal{S}(H)$ we have
$$B=\sup_{\nu \in \mathrm{co}\mathcal{S}(H)}\int_H U^{\mu}d\nu\leq \sup_{\nu\in\mathcal{R}(H)}\int_H U^{\mu}d\nu=A$$ $$\leq 3\varepsilon+\int_H f d\sigma\leq 3\varepsilon+\int_H U^{\mu} d\sigma \leq 3\varepsilon+\sup_{\sigma\in \mathrm{co}\mathcal{S}(H)}\int_H U^{\mu} d\sigma =3\varepsilon+B.$$
Since $\varepsilon$ was arbitrary, (2) is proved. By this computation and by Lemma 3 we have
$$\inf_{\mu \in \tilde{\mathcal{R}}(L)}\sup_{\nu\in \mathrm{Ex}\mathcal{R}(H)}E(\mu,\nu)=\inf_{\mu \in \tilde{\mathcal{R}}(L)}\sup_{\nu\in \mathrm{Ex}\mathcal{R}(H)}\int_HU^{\mu}d\nu$$ $$=\inf_{\mu \in \tilde{\mathcal{R}}(L)}\sup_{\nu\in \mathcal{R}(H)}\int_HU^{\mu}d\nu=\sup_{\nu\in \mathcal{R}(H)}\inf_{\mu \in \tilde{\mathcal{R}}(L)}E(\mu,\nu),$$
and together with Lemma 2 this ensures that
$$\sup_{\nu\in \mathcal{R}(H)}\inf_{\mu \in \mathrm{Ex}\tilde{\mathcal{R}}(L)}E(\mu,\nu)\leq \inf_{\mu \in \tilde{\mathcal{R}}(L)}\sup_{\nu\in \mathrm{Ex}\mathcal{R}(H)}E(\mu,\nu)= q\left(\tilde{\mathcal{R}}(L),\mathrm{Ex}\mathcal{R}(H)\right)$$ $$=\sup_{\nu\in \mathcal{R}(H)}\inf_{\mu \in \tilde{\mathcal{R}}(L)}E(\mu,\nu)\leq \sup_{\nu\in \mathcal{R}(H)}\inf_{\mu \in \mathrm{Ex}\tilde{\mathcal{R}}(L)}E(\mu,\nu)=\underline{q}\left(\mathcal{R}(H),\mathrm{Ex}\tilde{\mathcal{R}}(L)\right).$$

\medskip

\noindent {\bf Examples.}
Let $K\subset\subset X$, $\mathcal{R}(K)=\mathcal{M}_1(K)$ The extremal points of $\mathcal{R}(K)$ are the Dirac measures concentrated on the points of $K$. This example gives back the results of \cite{fare}. For further examples see the "Appendix".

\medskip

\lemma
Let $H\subset X, L\subset Y$ arbitrary and $k$ is a positive l.s.c. symmetric kernel. If for all $\mu \in \mathcal{R}(H)$ and $\varepsilon > 0$ there is a compact set $K(\varepsilon) \in \mathcal{K}(H)$ and a measure $\mu_{K(\varepsilon)}$ such that $\mu_{K(\varepsilon)} \in \mathcal{R}(K(\varepsilon))$ and $\mu_{K(\varepsilon)}\leq (1+\varepsilon)\mu|_{K(\varepsilon)}$, then
$$q(\mathcal{R}(H),\mathcal{S}(L))=\inf_{K\subset\subset H}q(\mathcal{R}(K),\mathcal{S}(L)).$$

\proof Obviously
$$q(\mathcal{R}(H),\mathcal{S}(L))=\inf_{\mu \in \mathcal{R}(H)}\sup_{\nu\in \mathcal{S}(L)}E(\mu,\nu)$$ $$\leq \inf_{\mu \in \mathcal{R}(K)\atop K\subset\subset H }\sup_{\nu\in \mathcal{S}(L)}E(\mu,\nu)=\inf_{K\subset\subset H }q(\mathcal{R}(K),\mathcal{S}(L)).$$
On the other hand let $\mu \in \mathcal{R}(H)$ and $\varepsilon > 0$ arbitrary.
$$\inf_{K\subset\subset H }q(\mathcal{R}(K),\mathcal{S}(L))\leq q(\mathcal{R}(K(\varepsilon)),\mathcal{S}(L))$$ $$\leq \sup_{\nu\in \mathcal{S}(L)}E(\mu_{K(\varepsilon)},\nu)\leq (1+\varepsilon)\sup_{\nu\in \mathcal{S}(L)}\int_L\int_{K(\varepsilon)}k(x,y)d\mu(x)d\nu(y)$$ $$\leq (1+\varepsilon)\sup_{\nu\in \mathcal{S}(L)}\int_L\int_{H}k(x,y)d\mu(x)d\nu(y).$$
Taking infimum over $\mu$ and tending to zero with $\varepsilon$ we obtain the converse inequality.

\medskip

Now we give some examples to energy-type expressions and investigate the assumptions of Lemma 4.

\noindent {\bf Examples.}

\noindent {\bf (1)} Let $\mathcal{R}(H)=\mathcal{M}_1(H)$, $\mathcal{R}(K)=\mathcal{M}_1(K)$, $\mu\in \mathcal{R}(H)$ arbitrary. Since $\mu$ is regular for all $\varepsilon > 0$ there is a $K(\varepsilon) \in \mathcal{K}(H)$ such that $\mu(K(\varepsilon))\geq \frac{1}{1+\varepsilon}$. Then $\mu_{K(\varepsilon)}:=\frac{\mu|_{K(\varepsilon)}}{\mu(K(\varepsilon))}$ fulfils the assumptions and we get back \cite[Lemma 2.4]{fare}.

\noindent {\bf (2)} First we remark that in this lemma the properties (e.g. the local compactness) of $Y$ does not play any role. We can apply this lemma to the Adams-Hedberg $p$-capacity of a set ($C_p$). Let $X=\mathbb{R}^n$, $L=Y=(M,\nu)$ a measure space. $1< p <\infty$. $\mathcal{S}(L):=\{\nu_f : d\nu_f=fd\nu, \ws f\geq 0, \|f\|_{\nu,p}\leq 1\}$. $K\subset \mathbb{R}^n$ compact. Then (cf. \cite[Theorem 2.5.1]{ah})
$$C_p(K)^{-\frac{1}{p}}=\min_{\mu\in \mathcal{M}_1(K)}\sup_{\nu_f \in \mathcal{S}(L)}E(\mu,\nu_f),$$
that is similarly to the previous example we have $C_p(H)^{-\frac{1}{p}}=\inf_{K\subset\subset H}C_p(K)^{-\frac{1}{p}}$.

\noindent {\bf (3)} The modified $p$-capacity cf. \cite[Definition 3.]{ho1}: $X$ is a locally compact Hausdorff space, $k$ is l.s.c. symmetric and positive on $X\times X$, $\mathcal{S}(X):=\{\nu_f : d\nu_f=fd\nu, \ws f\geq 0, \|f\|_{\nu,p}\leq 1\}$ $H\subset X$. The definition in the cited work uses probability measures which are compactly supported on $H$. Now we give the definition by measures concentrated on $H$, instead. Let
$$W_{p,0}(H)=\inf_{\mu\in \mathcal{M}_1(H)}\sup_{\nu\in\mathcal{S}(X)}E(\mu,\nu)$$
and
$$W_{p,\lambda}(H)=\inf_{\mu\in \mathcal{M}_1(H)}\sup_{\nu\in\mathcal{S}(X)}\left((1-\lambda)E(\mu,\nu)+\lambda E(\mu,\mu)\right). $$
Notice if $\lambda\neq 0$ then $W_{p,\lambda}(H)\neq q(\mathcal{M}_1(H),\mathcal{S}(X))$. Indeed $$W_{p,\lambda}(H)=\inf_{\mu\in \mathcal{M}_1(H)}\sup_{\nu\in\mathcal{S}(X)}E_J(\mu,\mu),$$ where $J(x,y)=(1-\lambda)\int_Xk(x,y)d\nu_f(y)+\lambda k(x,y)$, that is in this case we can not apply our discussion.

\noindent {\bf (4)} Let $H\subset X$, $f \geq 0$, a bounded and continuous function on $H$, and $c>0$ is a real number.
$$\mathcal{R}(H):=\left\{\mu \in \mathcal{M}_1(H) : \int_Hfd\mu\geq c\right\}.$$
(Assume, that $f$ and $c$ are given such that $\mathcal{R}(H)\neq\emptyset$.) Let $\mu \in \mathcal{R}(H)$, $\varepsilon>0$ arbitrary. We show that there is a compact set $K(\varepsilon)$ and a measure $\mu_{K(\varepsilon)}\in \mathcal{R}(K(\varepsilon))=\left\{\mu \in \mathcal{M}_1(K(\varepsilon)) : \int_Hfd\mu\geq c\right\}$ such that $\mu_{K(\varepsilon)}\leq (1+\varepsilon)\mu|_{K(\varepsilon)}$.

If $\mu$ is compactly supported on $H$ or $f\geq c$ $\mu$-a.e., there is nothing to prove. If there is a $0<\delta\leq \varepsilon$ and a compact set $K_{\delta}\subset H$ such that $\frac{1}{1+\varepsilon}\leq \mu(K_{\delta})\leq \frac{1}{1+\delta}$ and $f|_{K_{\delta}}\geq\frac{c}{1+\delta}$, then denoting by $K(\varepsilon)=K_{\delta}$, $\mu_{\varepsilon}:=\frac{\mu|_{K_{\delta}}}{\mu(K_{\delta})}$, $\mu_{\varepsilon}\leq (1+\varepsilon)\mu|_{K_{\delta}}$, and $\mu_{\varepsilon}\in \mathcal{R}(K(\varepsilon))$, since
$\int_{K_{\delta}}fd\mu_{\varepsilon}=\frac{1}{\mu(K_{\delta})}\int_{K_{\delta}}fd\mu\geq\frac{c}{(1+\delta)\mu(K_{\delta})}\geq c.$

Let $0<\delta\leq \varepsilon$ and $K_{\delta}\subset H$ compact such that $f|_{K_{\delta}}\geq\frac{c}{1+\delta}$. Assume that $1\geq \mu(K_{\delta})>\frac{1}{1+\delta}$ or $\mu(K_{\delta})<\frac{1}{1+\varepsilon}$. If there is a sequence $\delta_n \to 0$, such that $f|_{K_{\delta_n}}\geq\frac{c}{1+\delta_n}$ and $ \mu(K_{\delta_n})>\frac{1}{1+\delta_n}$, then $f\geq c$ $\mu$- a.e. on $H$, and by the regularity of $\mu$ we can choose $K(\varepsilon)$ as in the first example. So we can assume that there is a $0<\delta_0\leq \varepsilon$ so that if $0< \delta\leq \delta_0$ and there is a compact set $K_{\delta}$ such that $f|_{K_{\delta}}\geq\frac{c}{1+\delta}$, then $\mu(K_{\delta})\leq\frac{1}{1+\delta}$; moreover we can assume that $\mu(K_{\delta})\leq\frac{1}{1+\varepsilon}$.

Let $0< \delta\leq \delta_0$ $A_{\delta}:=\left\{x\in H : f(x)< \frac{c}{1+\delta}\right\}.$
By the assumption on $\mu$, the $\mu(H\setminus A_{\delta})<1$. If there is an $F\subset \subset A_{\delta}$ such that $0<\tilde{\varepsilon}=
\mu(A_{\delta}\setminus F)$, where $0<\tilde{\varepsilon}<\frac{\varepsilon}{1+2\varepsilon} (<\mu(A_{\delta}))$. Let $J\subset\subset H\setminus A_{\delta}$ such that $\int_{H\setminus A_{\delta}\setminus J}fd\mu\leq c \tilde{\varepsilon}\frac{\delta}{1+\delta}, \ws \ws \mu(H\setminus A_{\delta}\setminus J)\leq \tilde{\varepsilon}\frac{\delta}{1+\delta}.$
Let $K(\varepsilon)=F\cup J$, and $\mu_{\varepsilon}=\frac{\mu|_{K(\varepsilon)}}{\mu(K(\varepsilon))}$. Then
$\mu(A_{\delta})-\tilde{\varepsilon}+1-\mu(A_{\delta})-\tilde{\varepsilon}\frac{\delta}{1+\delta}\leq \mu(K(\varepsilon))\leq \mu(A_{\delta})-\tilde{\varepsilon}+1-\mu(A_{\delta})$,
and
$\int_{H\setminus K(\varepsilon)}fd\mu=\int_{ A_{\delta}\setminus F}fd\mu+\int_{H\setminus A_{\delta}\setminus J}fd\mu\leq \tilde{\varepsilon}\frac{c}{1+\delta}+c\tilde{\varepsilon}\frac{\delta}{1+\delta}=c\tilde{\varepsilon}$.
Thus
$\int_{K(\varepsilon)}fd\mu_{\varepsilon}=\frac{1}{\mu(K(\varepsilon))}\int_{K(\varepsilon)}fd\mu\geq c+\frac{(1-\mu(K(\varepsilon)))c-\int_{H\setminus K(\varepsilon)}fd\mu}{\mu(K(\varepsilon))}\geq c.$
That is $\mu_{\varepsilon}\in \mathcal{R}(K(\varepsilon))$, and $\frac{1}{\mu(K(\varepsilon))}\leq \frac{1}{1-\tilde{\varepsilon}\frac{1+2\delta}{1+\delta}}\leq \frac{1+\varepsilon}{1+\varepsilon-\tilde{\varepsilon}(1+2\varepsilon)}\leq 1+\varepsilon$.
If $\forall F\subset \subset A_{\delta}$ $\tilde{\varepsilon}<\frac{\varepsilon}{1+2\varepsilon}$ ensures that $\tilde{\varepsilon}=0$, then $K(\varepsilon)=K_{\delta}\cup F$ and $\mu_{K(\varepsilon)}=\mu|_{K(\varepsilon)}$.

\medskip

\remark
By Lemma 1, if $\mathcal{S}(L)$ is compact in the vague topology, $\forall \mu \in \mathcal{M}(H)$ $\inf_{\nu \in \mathcal{S}(L)}E(\mu,\nu)=\min_{\nu \in \mathcal{S}(L)}E(\mu,\nu)$.

\medskip

\lemma
Let $H\subset X$ arbitrary  and  $\mathcal{S}(L)$ is compact in the vague topology. If for an arbitrary measure $\mu \in \mathcal{R}(H)$ there is a net of measures $\{\mu_K\}_{\mathcal{K}_{\mu}(H)}\subset \mathcal{R}(H)$ such that $\mathcal{K}_{\mu}(H)\subset \mathcal{K}(H)$ and $\lim_{K\in\mathcal{K}_{\mu}(H)}\mu_K=\mu$ in the vague topology, then
$$\underline{q}(\mathcal{R}(H),\mathcal{S}(L))=\sup_{K\subset\subset H}\underline{q}(\mathcal{R}(K),\mathcal{S}(L)).$$

\proof
$$\underline{q}(\mathcal{R}(H),\mathcal{S}(L))=\sup_{\mu \in \mathcal{R}(H)}\inf_{\nu\in \mathcal{S}(L)}E(\mu,\nu)$$ $$\geq \sup_{\mu \in \mathcal{R}(K)\atop K\subset\subset H}\inf_{\nu\in \mathcal{S}(L)}E(\mu,\nu)=\sup_{K\subset\subset H}\underline{q}(\mathcal{R}(K),\mathcal{S}(L)),$$
always. On the other hand let $\mu \in \mathcal{R}(H)$ arbitrary again. Considering the net given by the assumption we define the corresponding family of measures.
Let $\nu_K$ be such that $\inf_{\nu \in \mathcal{S}(L)}E(\mu_K,\nu)=E(\mu_K,\nu_K)$. As $\mathcal{S}(L)$ is vaguely compact, there is a subnet $\mathcal{N}(H)\subset \mathcal{K}_{\mu}(H)$ such that $\lim_{K\in\mathcal{N}(H)}\nu_K=\nu_0$ in the vague topology with some $\nu_0$.
Since $E(\mu,\nu)$ is l.s.c. on $\mathcal{M}_1(H)\times \mathcal{N}(L)$,
$$\liminf_{\mathcal{N}(H)}E(\mu_K,\nu_K)\geq E(\mu,\nu_0)\geq \inf_{\nu \in \mathcal{S}(L)}E(\mu,\nu).$$
Thus
$$\sup_{K\subset\subset H}\underline{q}(\mathcal{R}(K),\mathcal{S}(L))=\sup_{K\subset\subset H\atop \mu \in\mathcal{R}(K)}\inf_{\nu \in \mathcal{S}(L)}E(\mu,\nu) \geq \liminf_{\mathcal{N}(H)}\inf_{\nu \in \mathcal{S}(L)}E(\mu_K,\nu)$$ $$=\liminf_{\mathcal{N}(H)}E(\mu_K,\nu_K)\geq \inf_{\nu \in \mathcal{S}(L)}E(\mu,\nu)$$
for all $\mu \in \mathcal{R}(H)$, so taking supremum in $\mu \in \mathcal{R}(H)$ we obtain the converse inequality.

\medskip

\noindent {\bf Examples.}

\noindent  {\bf (1)} Let $\mathcal{R}(H)=\mathcal{M}_1(H)$. Since $\mu(H)=1$ there is a $K_0\subset H$ such that $\mu(K_0)>0$. Let
$$\mathcal{K}_{\mu}(H):=\left\{K\subset\subset H : K_0\subset K\right\}.$$
Since $\mu$ is regular there exists a subnet $\{\mu_K\}_{K\in \mathcal{K}_{\mu}^*(H)}$ of the net indexed by the sets of $\mathcal{K}_{\mu}(H)$ ($\mu_K=\frac{\mu|_K}{\mu(K)}$) such that $\lim_{\mathcal{K}_{\mu}^*(H)}\mu_K=\mu $ in the vague topology.

\noindent  {\bf (2)} If $H$ is compact, Example (4) after Lemma 4 is a proper example here as well.

\medskip
\remark
Using Lemma 4 or Lemma 5, in special case we get back the Fuglede's-type result of the classical potential theory (cf. \cite[Theorem 5.3]{fare}).

\medskip

We extend some of the results of Fuglede and apply to our cases, cf. \cite[Theorem 2.3, Lemma 2.3.2, Theorem 2.4]{fu}. First we deal with the obvious direction.
We turn from  the vague topology to the $w^*$-topology, since in the applications the latter is the more useful. By Lemmas 4 and 5 we can assume that $H$ and $L$ are compact sets and in this case the vague topology coincides with the $w^*$-topology.

\begin{lemma}
Let $\mathcal{R}(H)\subset \mathcal{M}(H)$ be $w^*$-compact, convex set of measures. Let us assume that the kernel function is l.s.c. and positive. Then
$$w(\mathcal{R}(H))\leq q\left(\mathcal{R}(H)\right).$$\end{lemma}

\proof
Let $\mu\in \mathcal{R}(H)$ arbitrary. By the lower semicontinuity of the kernel $E(\mu)= \sup_{0\leq f\leq k\atop f\in C_c(X\times Y)}\int_H\int_Hf d\mu d\mu$. Since $\mathcal{R}(H)$ is  $w^*$-compact, convex there is a net in $\mathrm{Ex}\mathcal{R}(H)$, $\sigma_{\alpha}=\sum_{i=1}^{n_{\alpha}}\Theta_i\nu_i$ , $\sum_{i=1}^{n_{\alpha}}\Theta_i=1$ such that $\sigma_{\alpha} \stackrel{*}{\to} \mu$. That is
$$E(\mu)=\lim_{\alpha}\sup_{0\leq f\leq k\atop f\in C_c(X\times Y)}\int_H\int_Hfd\mu d\sigma_{\alpha} \leq \liminf_{\alpha}\sum_{i=1}^{n_{\alpha}}\Theta_i\int_H\int_H k(x,y)d\mu(x)d\nu_i(y)$$ $$\leq \liminf_{\alpha}\sum_{i=1}^{n_{\alpha}}\Theta_i \sup_{\nu\in \mathrm{Ex}\mathcal{R}(H)}\int_H\int_H k(x,y)d\mu(x)d\nu(y)=\sup_{\nu\in \mathrm{Ex}\mathcal{R}(H)}\int_H\int_H k(x,y)d\mu(x)d\nu(y).$$
Taking infimum over $\mu$, the lemma is proved.

\medskip

\remark

\noindent {\bf (1)} Let  $H\subset\subset X$. Notice that by Lemma 1, if a set of measures $\mathcal{R}(H)\subset \mathcal{M}(H)$ is $w^*$-compact, then there is an equilibrium measure $\mu_{\mathcal{R}(H)} \in \mathcal{R}(H)$ such that $E(\mu_{\mathcal{R}(H)})=w(\mathcal{R}(H))$.

\noindent {\bf (2)} A property fulfils nearly everywhere (n.e.) on a set $H$, if denoting by $N$ the set of points of $H$ for which the property does not fulfil, $w(\mathcal{M}_1(N))=\infty$. If a property fulfils n.e. on $H$, and a measure $\mu$ has finite Wiener energy, i.e. $E(\mu)<\infty$, then this property holds $\mu$-a.e. in $H$, cf. \cite[pages 153, 160]{fu}.

\medskip

\begin{defi}
A property fulfils $\mathcal{R}(H)$-nearly everywhere ($\mathcal{R}(H)$-n.e.) on a set $H$, if denoting by $N$ the set of points of $H$ for which the property does not fulfil, $w(\mathcal{R}(N))=\infty$, where $\mathcal{R}(N)=\{\mu \in \mathcal{R}(H): \mu \ws \mbox{is concentrated on } N \}$.
\end{defi}

\medskip

We need a lemma similar to \cite[Lemma 2.3.2]{fu}. Notice, that the opposite direction does not fulfil.

\lemma
Let $\mu\in \mathcal{M}(X)$, $H\subset X$, $0\leq t \leq \infty$, $k$ be a symmetric kernel, $\mathcal{R}(H)\subset \mathcal{M}_1(H)$. If

$$E(\nu,\mu)\geq t \ws \ws \ws \forall \ws \nu \in \mathcal{R}(H),  \ws E(\nu)<\infty, \ws  \mathrm{supp}\nu\subset\subset H,$$ then

$$U^{\mu}(x)\geq t \ws\ws \mathcal{R}(H)-\mathrm{n.e.} \ws x \in H.$$

\proof
Oppositely let us assume that $w(\mathcal{R}(N))<\infty$, where $N:=\{x\in H: U^{\mu}(x)<t\}$. That is there is a $\nu\in \mathcal{R}(N)\subset \mathcal{M}_1(N)$ with finite energy and supported on $N$ (cf. \cite[Lemma 2.3.1]{fu}). Then $E(\mu,\nu)=\int_NU^{\mu}d\nu<t$ which proves the statement by contradiction.

\medskip

\begin{defi}
Using the notation of the lemma we say that a set of measures $\mathcal{R}(H)\subset \mathcal{M}_1(H)$ is appropriate, if for all $\mu \in \mathcal{R}(H)$ with finite energy such that $E(\nu,\mu)\geq t$ for all $\nu \in \mathcal{R}(H)$, $E(\nu)<\infty$, $\mathrm{supp}\nu\subset\subset H$, it fulfils that $U^{\mu}(x)\geq t$ $\mu$-a.e. on $H$.
\end{defi}

\medskip

\begin{theorem}
Let $H,L \subset \subset X$, $\mathcal{R}(H)\subset \mathcal{M}_1(H)$ be $w^*$-compact, convex and appropriate set of measures and let $\mathcal{S}(L)\subset \mathcal{M}_1(L)$. Let us assume that the l.s.c. positive symmetric kernel function satisfies the Frostman's maximum principle. Then
$$w(\mathcal{R}(H))\geq q\left(\mathcal{R}(H),\mathcal{S}(L)\right).$$\end{theorem}

\proof
If $w(\mathcal{R}(H))=\infty$, the inequality is obvious, so let us assume that $w(\mathcal{R}(H))$ is finite.
Since $\mathcal{R}(H)$ is $w^*$-compact, by Lemma 1 there is an equilibrium measure $\mu_{\mathcal{R}(H)}\in \mathcal{R}(H)$. Let $\sigma\in \mathcal{R}(H)$ such that $E(\sigma)<\infty$, $\mathrm{supp}\sigma \subset\subset H$. Since $\mathcal{R}(H)$ is convex $\mu_{\Theta}:=(1-\Theta)\mu_{\mathcal{R}(H)}+\Theta\sigma\in \mathcal{R}(H)$ ($\Theta\in [0,1]$) and $F(\Theta)=E(\mu_{\Theta})$ has a minimum in $\Theta=0$. Thus $0$ has a right neighborhood where $F'\geq 0$. Since
$$\frac{1}{2}F'(\Theta)=-(1-\Theta)E(\mu_{\mathcal{R}(H)})+(1-2\Theta)E(\mu_{\mathcal{R}(H)},\sigma)+\Theta E(\sigma)\geq 0,$$
when $\Theta$ tends to zero, by the finiteness of the energy of $\sigma$ we have
$$E(\mu_{\mathcal{R}(H)},\sigma)\geq E(\mu_{\mathcal{R}(H)})\ws \forall \sigma\in \mathcal{R}(H), \ws E(\sigma)<\infty, \ws \mathrm{supp}\sigma \subset\subset H.$$
By Lemma 7 with $t=E(\mu_{\mathcal{R}(H)}$,
$U^{\mu_{\mathcal{R}(H)}}(x)\geq w(\mathcal{R}(H))$ $\mathcal{R}(H)$-n.e. on $H$ and by the assumption it fulfils $\mu_{\mathcal{R}(H)}$-almost everywhere. Let us assume now that there is an $x_0\in \mathrm{supp}\mu_{\mathcal{R}(H)}$ such that $U^{\mu_{\mathcal{R}(H)}}(x_0)> w(\mathcal{R}(H))$. By the lower semicontinuity there is a $\delta>0$ and a neighborhood of $G$ $x_0$, such that for all $x\in G\cap \mathrm{supp}\mu_{\mathcal{R}(H)}$ $U^{\mu_{\mathcal{R}(H)}}(x_0)\geq w(\mathcal{R}(H))+\delta$. Then
$$w(\mathcal{R}(H))=\int_GU^{\mu_{\mathcal{R}(H)}}d\mu_{\mathcal{R}(H)}+\int_{H\setminus G}U^{\mu_{\mathcal{R}(H)}}d\mu_{\mathcal{R}(H)}$$ $$\geq (w(\mathcal{R}(H))+\delta)\mu_{\mathcal{R}(H)}(G)+w(\mathcal{R}(H))\mu_{\mathcal{R}(H)}(H\setminus G),$$
that is $\mu_{\mathcal{R}(H)}(G)=0$ so $x_0\notin \mathrm{supp}\mu_{\mathcal{R}(H)}$. Since the kernel fulfils the Frostman's maximum principle $U^{\mu_{\mathcal{R}(H)}}(x)\leq w(\mathcal{R}(H))$ for all $x\in H$. Because $\mathcal{S}(L)\subset \mathcal{M}_1(L)$, for every $\nu \in \mathcal{S}(L)$
$$\int_L U^{\mu_{\mathcal{R}(H)}}d\nu\leq w(\mathcal{R}(H)).$$
Taking supremum over $\mathcal{S}(L)$ and infimum over $\mathcal{R}(H)$ the inequality is proved.

\medskip

\begin{corollary}

Let $H\subset\subset X$, $\mathcal{R}(H)\subset \mathcal{M}_1(H)$ be appropriate, $w^*$-compact, convex set of measures and the l.s.c. positive symmetric kernel function $k$ satisfies the Frostman's maximum principle. Then
$$w(\mathcal{R}(H))= q\left(\mathcal{R}(H)\right).$$\end{corollary}

\medskip

It is clear that $\mathcal{M}_1(H)$ is appropriate, and in this case we get back the original theorem of Fuglede. Below we give an example for appropriate set of measures different from $\mathcal{M}_1(H)$.

\medskip

\noindent {\bf Example.}
Let $K\subset\subset X$, ($\mathrm{card}K=\infty$) $f$ is a real valued, positive, continuous and for simplicity totally non constant function on $K$, that is for all $G\subset K$ open set $f|_G \not\equiv \mathrm{const}$. The kernel $k$ is positive, symmetric and l.s.c. and infinite at the diagonal. Let
$$\mathcal{R}(K)=\{\mu\in \mathcal{M}_1(K):\int_Xfd\mu=c\}.$$
We choose $c$ such that $\mathcal{R}(K)\neq \emptyset$.
Let us assume that $\mu \in \mathcal{R}(K)$ with finite energy $w$, $0\leq w<\infty$, $E(\nu,\mu)\geq w$ for all $\nu \in \mathcal{R}(K)$, $E(\nu)<\infty$, $\mathrm{supp}\nu\subset K$. Let $N:=\{x\in \mathrm{supp}\mu : U^{\mu}(x)<w\}$. Consider $\mathcal{R}(N)=\{\mu\in \mathcal{M}_1(N):\int_Xfd\mu=c\}$. By Lemma 7 $w(\mathcal{R}(N))=\infty$. Thus $\mu(N)\neq 1$. We show that $\mu(N)=0$. Let us assume the $1>\mu(N)>0$, that is $w(\mathcal{M}_1(N))<\infty=w(\mathcal{R}(N))$. It follows that for all $\nu\in \mathcal{M}_1(N)$, $E(\nu)<\infty$ $\int_Nfd\nu>c$ or for all $\nu\in \mathcal{M}_1(N)$ $E(\nu)<\infty$ $\int_Nfd\nu<c$, since if there are two measures $\nu_1$, $\nu_2$ with finite energy for which the integrals have opposite sign, their convex combination gives a measure with finite energy in $\mathcal{R}(N)$. So assume that for all $\nu\in \mathcal{M}_1(N)$, $E(\nu)<\infty$ $\int_Nfd\nu>c$, say.  Except a set of zero capacity, that is a except set $Z$ with $w(\mathcal{M}_1(Z))=0$, we can divide $\mathrm{supp}\mu$ into four subsets: $S_1:=\{x\in \mathrm{supp}\mu : U^{\mu}(x)< w, \ws f>c\}$, $S_2:=\{x\in \mathrm{supp}\mu : U^{\mu}(x)\geq w, \ws f<c\}$, $S_3:=\{x\in \mathrm{supp}\mu : U^{\mu}(x)\geq w, \ws f>c\}$, $S_4:=\{x\in \mathrm{supp}\mu : U^{\mu}(x)< w, \ws f<c\}$. ($\mathrm{cap}S_1>0$.) If $\mu(S_i)\neq 0$, let $\mu_i:=\frac{\mu|_{S_i}}{\mu(S_i)}$ $i=1,2,3,4$.
If there is a $\sigma \in \mathcal{M}_1(S_4)$ with finite energy, there is an $\alpha \in (0,1)$ such that $\alpha\int_{S_1}fd\mu_1 +(1-\alpha)\int_{S_4}fd\sigma=c$ and so $\nu=\alpha \mu_1 +(1-\alpha)\sigma \in \mathcal{R}(K)$ and $E(\mu,\nu)<w$. That is $\mathrm{cap}S_4=0$. Assume that $\mu(S_i)>0$ $i=1,2,3$.  Let us define $p_i, q_1> 0$, $i=1,2,3$, $q_2, q_3\geq 0$ as follows. $\int_{S_1}fd\mu_1=c+p_1$, $\int_{S_1}U^{\mu}d\mu_1=w-q_1$, $\int_{S_2}fd\mu_2=c-p_2$, $\int_{S_2}U^{\mu}d\mu_2=w+q_2$, $\int_{S_3}fd\mu_3=c+p_3$, $\int_{S_3}U^{\mu}d\mu_3=w+q_3$. Since $\mu \in \mathcal{R}(K)$, there are $\alpha, \beta, \gamma > 0$  solutions of the system of linear equations:
$$\left[\begin{array}{ccc} 1& 1&1 \\ p_1 & -p_2&p_3\\-q_1&q_2&q_3\end{array}\right]\left[\begin{array}{c}\alpha\\ \beta\\ \gamma\end{array}\right]=\left[\begin{array}{c}1\\ 0\\r\end{array}\right]$$
with $r=0$. A positive solution of the system above with $r<0$ means that there is a positive measure $\nu \in \mathcal{R}(K)$ such that $E(\mu,\nu)<w$. Let $D_1=\left|\begin{array}{cc} p_1&-p_2\\-q_1&q_2\end{array}\right|$, $D_2=\left|\begin{array}{cc} -p_2&p_3\\q_2&q_3\end{array}\right|$, $D_3=\left|\begin{array}{cc} p_1&p_3\\-q_1&q_3\end{array}\right|$. Suppose, $D:=D_2-D_3+D_1\neq 0$. Since $D_3>0$, $-D_2\geq 0$, the solutions $\alpha=\frac{D_2}{D}$, $\beta =\frac{-D_3}{D}$, $\gamma=\frac{D_1}{D}$ cannot be positive if $D>0$. The solutions are positive, if $D<0$ and $D_1, D_2< 0$. In this case with an $\varepsilon>0$ small enough, we have a positive solution of the linear system above with $r=-\varepsilon$, $\alpha=\frac{D_2+\varepsilon(p_3+p_2)}{D}$, $\beta =\frac{-D_3+\varepsilon(p_1-p_3)}{D}$, $\gamma=\frac{D_1-\varepsilon(p_1+p_2)}{D}$, it leads to a contradiction. If $D=0$, we have a nonnegative solution only if $D_1=D_2=D_3=0$, which is a contradiction. So oppositely to the assumption, at least one of the three sets has zero $\mu$-measure. If $\mu(S_1)=0$, then $\mu(\{x\in \mathrm{supp}\mu : U^{\mu}(x)>w\})=0$ and so $\mu(N)=0$ which is a contradiction. By the assumption on $f$, $\mu(S_2)$ can not be zero. That is $\mu$ is concentrated on $S_1\cup S_2$. So $\left[\begin{array}{cc} p_1 & -p_2\\-q_1&q_2\end{array}\right]\left[\begin{array}{c}\alpha\\ 1- \alpha\\ \end{array}\right]=\left[\begin{array}{c} 0\\r\end{array}\right]$ with $r=0$ has a solution $\alpha \in (0,1)$ that is $D_1=0$. Since $r=\frac{D_1}{p_1+p_2}$, the sign of $D_1$ and $r$ are coincide. We want to reduce $D_1$ by some perturbation of $\mu_i$-s. Let $h_i: S_i \to (0,2)$, such that $\int_{S_i}h_id\mu_i=1$. Let $d\mu_{1,\Theta}=((2\Theta-1) h_1+2(1-\Theta))d\mu_1$, $d\mu_{2,\Psi}=((2\Psi-1)h_2+2(1-\Psi))d\mu_2$, $\int_{S_1}fd\mu_{1,1}=c+a$; $\int_{S_1}U^{\mu}d\mu_{1,1}=w-d$; $\int_{S_2}fd\mu_{2,1}=c-u$; $\int_{S_2}U^{\mu}d\mu_{2,1}=w+x$, $(a,d,u>0)$, $(x\geq 0)$. We can assume that $\left|\begin{array}{cc} p_1 & -u\\-q_1&x\end{array}\right|\neq 0$, say. Then the equation is $\left[\begin{array}{cc}(2\Theta-1)a+2(1-\Theta) p_1 & -((2\Psi-1)u+2(1-\Psi)p_2)\\-((2\Theta-1)d+2(1-\Theta)q_1)&(2\Psi-1)x+2(1-\Psi)q_2\end{array}\right]\left[\begin{array}{c}\alpha\\ 1- \alpha\\ \end{array}\right]=\left[\begin{array}{c} 0\\r\end{array}\right]$.
The determinant, $D_1(\Theta,\Psi)=0$ if $\Theta=\Psi=\frac{1}{2}$, moreover $\forall \Theta, \Psi, \alpha \in(0,1)$ $\alpha \mu_{1,\Theta}+(1-\alpha)\mu_{2,\Psi}\in \mathcal{M}_1(K)$. On the other hand $D_1(\Theta,\Psi)=4\Theta\Psi((a-p_1)(x-q_2)-(d-q_1)(u-p_2))-2\Theta((a-p_1)(x-2q_2)-(d-q_1)(u-2p_2))-2\Psi((a-2p_1)(x-q_2)-(d-2q_1)(u-p_2))+((a-2p_1)(x-2q_2)-(d-2q_1)(u-2p_2))$. It is clear that this expression does not have a local minimum; if the coefficient of $\Theta\Psi$ is nonzero, the expression $D_1(\Theta,\Psi)=z$ is a hyperbolic paraboloid, if the coefficient of $\Theta\Psi$ is zero and the linear part is nonzero, then it is a skew plane, the point $\left(\frac{1}{2},\frac{1}{2}\right)$ is in the level $z=0$, so we can move the variables a little to be $D_1(\Theta,\Psi)<0$ and $\alpha(\Theta,\Psi)>0$, which is a contradiction. If all the coefficients are zero, the constant term has to be zero as well, but it is impossible by the assumption on $h_1$. Thus oppositely to our original assumption, $\mu(N)=0$.

\subsection{Discretization}

In the classical potential theory the $n$th Chebyshev constant is a discrete version of the potential, and when the kernel satisfies the Frostman's maximum principle, it coincides with the Wiener-energy. Similar theorems were proved in more general cases for instance in \cite{fana}, \cite{fare}, \cite{fu}, \cite{ho1}. In this subsection we extend the above-mentioned results to this more general structure. Under certain conditions it supports the discretization of the infinite quadratic programming problems of the next section.

\begin{defi}
Let $H \subset X$, $L \subset Y$. The $n$th Chebyshev constant of $\mathcal{R}(H)$ related to $\mathcal{S}(L)$ is
$$\underline{q}\left(\mathcal{R}_{n}(H),\mathcal{S}(L)\right),$$
where
$$\mathcal{R}_{n}(H):=\{\mu \in\mathcal{R}(H): \mathrm{card}(\mathrm{supp}\mu)\leq n\}.$$
The $n$th dual Chebyshev constant is
$$q\left(\mathcal{R}_{n}(H),\mathcal{S}(L)\right).$$\end{defi}

\medskip

Let us observe, that the $n$th Chebyshev constant and the $n$th dual Chebyshev constant have limit when $n$ tends to infinity.
For all $\Theta \in (0,1)$, $\mu_m \in \mathcal{R}_{m}(H)$, $\mu_n\in \mathcal{R}_{n}(H)$ we have
$$\underline{q}(\mathcal{R}_{m+n}(H),\mathcal{S}(L))=\sup_{\mu_{m+n}\in \mathcal{R}_{m+n}(H)}\inf_{\nu \in\mathcal{S}(L)}E(\mu,\nu)$$ $$\geq \sup_{\Theta\mu_{m}+(1-\Theta)\mu_n \atop {\Theta \in (0,1)\atop \mu_m\in \mathcal{R}_{m}(H),\mu_n\in \mathcal{R}_{n}(H)}}\inf_{\nu \in\mathcal{S}(L)}E(\Theta \mu_{m}+(1-\Theta)\mu_n,\nu)\geq \inf_{\nu \in\mathcal{S}(L)}E(\Theta \mu_{m}+(1-\Theta)\mu_n,\nu)$$ $$\geq \Theta\inf_{\nu \in\mathcal{S}(L)}E(\mu_{m},\nu)+(1-\Theta)\inf_{\nu \in\mathcal{S}(L)}E(\mu_n,\nu).$$
Let us take supremum over $\mathcal{R}_{m}(H)$ and $\mathcal{R}_{n}(H)$ respectively. Choosing $\Theta=\frac{m}{m+n}$ we get that $\underline{q}(\mathcal{R}_{n}(H),\mathcal{S}(L))$ is quasi-increasing, that is it has a limit.

Similarly $q(\mathcal{R}_{n}(H),\mathcal{S}(L))$ is quasi-decreasing that is it also has a limit.

\medskip

\begin{defi}
The Chebyshev constant of $\mathcal{R}(H)$ related to $\mathcal{S}(L)$ is
$$M\left(\mathcal{R}(H),\mathcal{S}(L)\right):=\lim_{n\to \infty}\underline{q}\left(\mathcal{R}_{n}(H),\mathcal{S}(L)\right),$$
and the dual Chebyshev constant is
$$\overline{M}\left(\mathcal{R}(H),\mathcal{S}(L)\right):=\lim_{n\to \infty}q\left(\mathcal{R}_{n}(H),\mathcal{S}(L)\right).$$
As above if $\mathcal{R}(H)$ is $w^*$-compact, convex, let us denote by
$$M\left(\mathcal{R}(H)\right):= M\left(\mathcal{R}(H),\mathrm{Ex}\mathcal{R}(H)\right), \ws \ws\overline{M}\left(\mathcal{R}(H)\right):= \overline{M}\left(\mathcal{R}(H),\mathrm{Ex}\mathcal{R}(H)\right).$$
\end{defi}

\medskip

\theorem
Let $H\subset\subset X$, $L\subset\subset Y$ such that $\mathcal{R}(H)$ is $w^*$-compact, convex and $\mathcal{S}(L)$ is $w^*$-compact. Then
$$M\left(\mathcal{R}(H),\mathcal{S}(L)\right)=\underline{q}\left(\mathcal{R}(H),\mathcal{S}(L)\right).$$

\proof
For every $n$, $\underline{q}\left(\mathcal{R}_{n}(H),\mathcal{S}(L)\right)\leq \underline{q}\left(\mathcal{R}(H),\mathcal{S}(L)\right)$ obviously, thence \\$M\left(\mathcal{R}(H),\mathcal{S}(L)\right)\leq\underline{q}\left(\mathcal{R}(H),\mathcal{S}(L)\right)$.

On the other hand since $\mathcal{R}(H)$ is $w^*$-compact and convex, for all $ \mu\in \mathcal{R}(H)$ there is a net of measures $\{\mu_{\alpha}\}_A\subset \mathrm{co Ex}\mathcal{R}(H)$ such that $\mu_{\alpha} \stackrel{*}{\to} \mu$ and $\mu_{\alpha}=\sum_{i=1}^{m_{\alpha}}\Theta_i\lambda_i$, where $\sum_{i=1}^{m_{\alpha}}\Theta_i=1$ and $\lambda_i\in \mathrm{Ex}\mathcal{R}(H)$, $i=1,\dots ,m_{\alpha}$. As $\mathcal{S}(L)$ is $w^*$-compact it can be proved similarly to Lemma 5, that for each $\mu_{\alpha}$ there exists a $\nu_{\alpha}\in \mathcal{S}(L)$, such that
\begin{equation}\int_L\int_Hk(x,y)d\mu_{\alpha}(x)d\nu_{\alpha}(y)=\inf_{\sigma\in \mathcal{S}(L)}\int_L\int_Hk(x,y)d\mu_{\alpha}(x)d\sigma(y).\end{equation}
Again by the $w^*$-compactness of $\mathcal{S}(L)$ there is a $w^*$-convergent subnet of $\{\nu_{\alpha}\}$ (denoted by $\{\nu_{\alpha}\}$ again) which tends to $\nu$. According to the lower semicontinuity of the function $(\mu,\nu)\to E(\mu,\nu)$ (cf. Lemma 1)
\begin{equation}\int_L\int_H k(x,y)d\mu(x)d\nu(y)\leq \liminf_{\alpha} \int_L\int_Hk(x,y)d\mu_{\alpha}(x)d\nu_{\alpha}(y).\end{equation}
From (7) and (8) we have
$$\inf_{\sigma \in \mathcal{S}(L)}\int_L\int_H k(x,y)d\mu(x)d\sigma(y)\leq \int_L\int_H k(x,y)d\mu(x)d\nu(y)$$ $$\leq \liminf_{\alpha} \int_L\int_Hk(x,y)d\mu_{\alpha}(x)d\nu_{\alpha}(y)=\liminf_{\alpha}\inf_{\sigma \in \mathcal{S}(L)}\int_L\int_H k(x,y)d\mu_{\alpha}(x)d\sigma(y)$$ $$\leq \liminf_{\alpha}\sup_{\lambda \in \mathcal{R}_{m_{\alpha}}(H)}\inf_{\sigma \in \mathcal{S}(L)}\int_L\int_H k(x,y)d\lambda(x)d\sigma(y)$$ $$=\liminf_{\alpha}\underline{q}\left(\mathcal{R}_{m_{\alpha}}(H),\mathcal{S}(L)\right)=M\left(\mathcal{R}(H),\mathcal{S}(L)\right).$$
Taking supremum over $\mathcal{R}(H)$ we have the opposite inequality.

\medskip

Summarizing, we arrived to the equality  
$$M\left(\mathcal{R}(H)\right)=\underline{q}\left(\mathcal{R}(H)\right)=q\left(\mathcal{R}(H)\right)=w\left(\mathcal{R}(H)\right),$$ 
provided that $H\subset\subset X$; $\mathcal{R}(H)$ is appropriate, $w^*$-compact, convex; and $k$ satisfies the Frostman's maximum principle. Now we are in position to apply these results to transform an infinite quadratic programming problem to a semi-infinite or to a finite one.

\section{Infinite Quadratic Programming}

In this section we study infinite dimensional quadratic programming problems by potential-theoretic methods. To discretize the problem we apply both cutting plane-and Chebyshev-type methods. In the Appendix we give several examples on sets of measures and kernels for which the results of the previous section can be applied.

The continuous/infinite quadratic programming problem was examined by several authors in metric spaces and in $L^p$ spaces on compact interals, cf. e.g. \cite{w}, \cite{cw}, \cite{wfl}. The problem in metric space, say is as it follows.
Let $X,Z$ be compact, metric spaces, $\Phi(x,z): X\times Z\to \mathbb{R}$, $g: Z\to \mathbb{R}$, $f(x,y): X\times X\to \mathbb{R}$, $h(x): X\to \mathbb{R}$ be continuous functions, furthermore let $f(x,y)$ be symmetric and positive semi-definite. Consider
\begin{equation}\inf_{\mu\in \mathcal{R}_{\Phi,g}^+(X)}\frac{1}{2}\int_X\int_X f(x,y) d\mu(x)d\mu(y)+\int_Xh(x)d\mu(x),\end{equation}
where
\begin{equation}\mathcal{R}_{\Phi,g}^+(X):=\{\mu\in \mathcal{M}(X): \int_X\Phi(x,z)d\mu(x)\geq g(z) \ws \ws \forall z\in Z\}.\end{equation}
It is clear that if  the total variation of the the optimal measure is between two positive constants $c_1$ and $c_2$, say then it is enough to deal with measures for which  $c_1\leq \mu(X)\leq c_2$. Factoring out $\mu(X)=c$ we can reformulate the problem above as it follows.
\begin{equation}\inf_{c_1\leq c\leq c_2}\inf_{\mu\in \mathcal{R}_{\Phi,g}(X)}c^2\int_X\int_X k_c(x,y)d\mu(x)d\mu(y)-c^2K_c,\end{equation}
where
\begin{equation}\mathcal{R}_{\Phi,g}(X)^c:=\{\mu\in \mathcal{M}_1(X): \int_Xc\Phi(x,z)d\mu(x)\geq g(z) \ws \ws \forall z\in Z\},\end{equation}
and the symmetric continuous positive kernel function is
\begin{equation}k_c(x,y)=\frac{1}{2}\left(f(x,y)+\frac{1}{c}(h(x)+h(y))\right)+K_c\end{equation}
with a suitable positive constant $K_c$.

Hereinafter we deal with the second minimum, and we omit the notation of dependence on $c$. That is we are interested in the following minimum problems.

Let $X,Z$ be compact, metrizable spaces, $\Phi(x,z): X\times Z\to \mathbb{R}$, $g: Z\to \mathbb{R}$ be continuous functions, $k(x,y)$ be a positive symmetric lower semicontinuous kernel on $X\times X$.
$$\mathrm{ (CQP_i)}:\ws\ws\ws \inf_{\mu\in \mathcal{R}_{\Phi,g}(X)}\int_X\int_X k(x,y)d\mu(x)d\mu(y)=w\left(\mathcal{R}_{\Phi,g}(X)\right), \ws\ws \mbox{where}$$
\begin{equation}  \mathcal{R}_{\Phi,g}(X):=\{\mu\in \mathcal{M}_1(X): \int_X\Phi(x,z)d\mu(x)\geq g(z) \ws \ws \forall z\in Z\}.\end{equation}
$$\mathrm{ (CQP_e)}:\ws\ws\ws \inf_{\mu\in \mathcal{R}_{\Phi,g}^e(X)}\int_X\int_X k(x,y)d\mu(x)d\mu(y)=w\left(\mathcal{R}_{\Phi,g}^e(X)\right), \ws\ws \mbox{where}$$
\begin{equation} \mathcal{R}_{\Phi,g}^e(X):=\{\mu\in \mathcal{M}_1(X): \int_X\Phi(x,z)d\mu(x)= g(z) \ws \ws \forall z\in Z\}.\end{equation}

Following the logic of E. Anderson (cf. \cite{a}), we examine mainly $\mathrm{ (CQP_e)}$. Our other motivation to deal with the "equality problem" is that there are nice methods to discuss $\mathrm{ (CQP_i)}$ under Slater condition, by duality (cf. \cite{wwt}).

One of our main tools is the cutting plane algorithm.

\subsection{The cutting plane or greedy algorithm}

We follow the chain of ideas of \cite[Section 2]{cw}. Here we discuss "equality" and "inequality " cases. We deal with the "equality case". The algorithm in the "inequality case" can be found basically in \cite{cw} and similar to the determination of the so-called greedy energy or Leja points (cf. \cite{ho1} and the references therein).

Let us assume that $\mathcal{R}_{\Phi,g}^e(X)\neq \emptyset$.
Let $n=1$, choose any $z_1, z_2 \in Z$, and $Z_2=\{z_1,z_2\}$. Next let $Z_{2n}=\{z_1,z_2,\dots ,z_{2n}\}$ be already defined, let $\mu_n$ be the "minimal measure" with respect to $\mathcal{R}_{\Phi,g,Z_{2n}}^e(X):=\{\mu\in \mathcal{M}_1(X): \int_X\Phi(x,z)d\mu(x)= g(z) \ws z\in Z_{2n}\}$, that is
$$E(\mu_n)= \inf_{\mu\in \mathcal{R}_{\Phi,g,Z_{2n}}^e(X)}E(\mu).$$
Observe that $\mathcal{R}_{\Phi,g,Z_{2n}}^e(X)=\{\mu\in \mathcal{M}_1(X): \int_X\Phi(x,z_i)d\mu(x) = g(z_i)\ws i=1, \dots, 2n\} \\= \{\mu\in \mathcal{M}_1(X): \int_Xf_i(x)d\mu(x)= c_i \ws i=1, \dots, 2n\}$, where $f_i(x)=\Phi(x,z_i)$ and $c_i=g(z_i)$. Since $f_i$ are continuous, $\mathcal{R}_{\Phi,g,Z_{2n}}^e(X)$ is a $w^*$-closed subset of a $w^*$-compact set, thus the infimum is a minimum. The situation is the same in the "inequality case".
Let us denote by
$$\Psi_n(z)=\int_X\Phi(x,z)d\mu_n(x)- g(z).$$
If $\Psi_n(z)\equiv 0$ then stop. Otherwise we continue the procedure with determination of $z_{2n+1}$ and $z_{2n+2}$. Since $\Psi_n$ is continuous and $Z$ is compact, there exist $z_{2n+1}, z_{2n+2} \in Z$ such that
$$\Psi_{n}(z_{2n+1})=\min_{z\in Z}\Psi_{n}(z), \ws\ws \Psi_{n}(z_{2n+2})=\max_{z\in Z}\Psi_{n}(z).$$
Then $Z_{2n+2}=Z_{2n}\cup \{z_{2n+1},z_{2n+2}\}$. In the "inequality case", if $\Psi_n(z)\geq 0$ then stop, and otherwise we follow the "$n$ is odd" case (cf. \cite{cw}).

\medskip

\begin{lemma}
Supposing that  $\mathcal{R}_{\Phi,g}(X)\neq \emptyset$ $\left(\mathcal{R}_{\Phi,g}^e(X)\neq \emptyset\right)$, the cutting plane algorithm has a limit, that is
$$\lim_{n\to\infty}w\left(\mathcal{R}_{\Phi,g,Z_{2n}}^e(X)\right)=w\left(\mathcal{R}_{\Phi,g}^e(X)\right).$$
\end{lemma}

\medskip

\noindent To prove the lemma we need the following result.

\medskip

\lemma\cite[Lemma 3.11]{fare}
If $X$ is compact, $\Phi(x,z)$ is continuous on $X\times Z$, the mapping
$$M : \mathcal{M}_1(X)\to C(X), \ws\ws \mu \mapsto U^{\mu}_{\Phi}=\int_X\Phi(\cdot,z)d\mu(x)$$
is continuous from the $w^*$-topology to the sup-norm topology.

\medskip

\proof (of Lemma 8) The "equality case".
If the algorithm is finite, there is nothing to prove. We deal with the infinite case.
Since $\mathcal{R}_{\Phi,g}^e(X)\subset \dots \subset \mathcal{R}_{\Phi,g,Z_{2n+2}}^e(X)\subset \mathcal{R}_{\Phi,g,Z_{2n}}^e(X)\subset \dots \subset \mathcal{R}_{\Phi,g,Z_{2}}^e(X)$, thus $E(\mu_n)$ is increasing and for all $n$ $E(\mu_n)\leq w(\mathcal{R}_{\Phi,g}^e(X))$ so it has a limit in the extended sense. Assume that $w\left(\mathcal{R}_{\Phi,g}^e(X)\right)<\infty$. Let us suppose contrary that with some $\eta >0$ $\lim_{n\to\infty}E(\mu_n)=w(\mathcal{R}_{\Phi,g}^e(X))-\eta$. Since $\{\mu_n\}\subset \mathcal{M}_1(X)$, it has a $w^*$-convergent subsequence, $\mu_{n_k}\stackrel{*}{\to}\mu$ and $\mu \in \mathcal{M}_1(X)$ as well. Since  $k$  is l.s.c.
$$E(\mu)=\sup_{0\leq h(x,y)\leq k(x,y)\atop h\in C_c(X\times X)}\int_{X\times X}h(x,y)d\mu(x)\times\mu(y)$$ $$=\sup_{0\leq h(x,y)\leq k(x,y)\atop h\in C_c(X\times X)}\lim_{k\to\infty}\int_X\int_X h(x,y)d\mu_{n_k}(x)d\mu_{n_k}(y)$$ $$\leq \sup_{0\leq h(x,y)\leq k(x,y)\atop h\in C_c(X\times X)}\liminf_{k\to\infty}\int_X\int_X k(x,y)d\mu_{n_k}(x)d\mu_{n_k}(y)$$ $$=\lim_{n\to\infty}E(\mu_n)=w(\mathcal{R}_{\Phi,g}(X))-\eta.$$
We show that $\mu\in \mathcal{R}_{\Phi,g}^e(X)$, which leads to a contradiction. \\
Let $\Psi(z)=\int_X\Phi(x,z)d\mu(x)- g(z)$ and let us define $z^*, z^{**}\in Z$ such that $\Psi(z^*)=\min_{z\in Z}\Psi(z)$, $\Psi(z^{**})=\max_{z\in Z}\Psi(z)$. Because $Z$ is compact $\{n_k\}$ has a subsequence such that $\{z_{2n_{k_l}+2}\}\to z_e$ and a (possibly different) subsequence so that $\{z_{2n_{k_j}+1}\}\to z_o$. The construction ensures that $\Psi(z_i)= 0$ for all $i\geq 1$, thus by the continuity of $\Psi$, $\Psi(z_e)=\Psi(z_o)=  0$. Also by the construction $\Psi_{n}(z^*)\geq \Psi_{n}(z_{2n+1})$, $n\geq 1$ and $\Psi_{n}(z^{**})\leq \Psi_{n}(z_{2n+2})$. Thus
\begin{equation}\Psi_{n_{k_l}}(z^{**})\leq \Psi_{n_{k_l}}(z_{2n_{k_l}+2}), \ws\ws l=1,2,\dots,\end{equation}
and
\begin{equation}\Psi_{n_{k_j}}(z^*)\geq \Psi_{n_{k_j}}(z_{2n_{k_j}+1}), \ws\ws j=1,2,\dots.\end{equation}
Since $\Psi_{n_k}$ is obviously pointwise convergent, the left-hand side of (16) and (17) tend to $\Psi(z^{**})$ and $\Psi(z^*)$ respectively. On the other hand $|\Psi_{n_{k_l}}(z_{2n_{k_l}+2})- \Psi(z_e)|\leq |\Psi_{n_{k_l}}(z_{2n_{k_l}+2})- \Psi(z_{2n_{k_l}+2})|+|\Psi(z_{2n_{k_l}+2})- \Psi(z_e)|$. Thus according to Lemma 9, and by the continuity of $\Psi$, $\Psi_{n_{k_l}}(z_{2n_{k_l}+2})\to \Psi(z_e)$ and similarly $\Psi_{n_{k_j}}(z_{2n_{k_j}+1})\to \Psi(z_o)$. That is
$$0= \Psi(z_o)\leq \Psi(z^*), \ws\ws \mbox{and}\ws\ws 0= \Psi(z_e)\geq \Psi(z^{**}),$$
which proves the statement. The proof of the "inequality case" is the same following the odd indices.

\medskip

Our aim is to give some semi-infinite or finite methods by the Chebyshev and the cutting plane algorithm which solve $(CQP_e)$, in limit.

\medskip

First, following the logic of the example after Corollary 2, we show that the sets of measures $\mathcal{R}_{\Phi,g,Z_{2n}}^e(X)$ are appropriate.

\medskip

\begin{lemma}Let $f_1, \dots, f_n$ be positive, continuous, totally nonconstant functions on $X$, $k$ is a symmetric l.s.c. kernel, let  $c_i$ be positive numbers, $i=1, \dots, n$. Assume that $\mathcal{M}_{\mathcal{F}}(X):=\{\mu \in \mathcal{M}_1(X) : \int_Xf_id\mu =c_i \ws i=1, \dots, n\}\neq \emptyset$. If $k$ is infinite at the diagonal or $1,f_1, \dots, f_n$ is a strong Chebyshev (or Descartes) system (see e.g. \cite{eb}), then $\mathcal{M}_{\mathcal{F}}(X)$ is appropriate.\end{lemma}

\proof Assume that $w:=w(\mathcal{M}_{\mathcal{F}}(X))=E(\mu)<\infty$. Now by Theorem 2, for any $\sigma\in \mathcal{M}_{\mathcal{F}}(X)$ $E(\mu,\sigma)\geq E(\mu)$. Let us suppose indirectly that $\mu(\{x\in X:U^{\mu}(x)<w\})>0$. Let $S_1:=\{x\in X:U^{\mu}(x)<w\}$, say and (without taking any care of the method) let $S_1, S_2, \dots, S_{n+1}$ be a partition of $X$ such that $\mu(S_k)>0$ $k=1, \dots, n+1$. It can be done, since if $k$ is infinite at the diagonal, a measure with finite energy is not atomic, and if $f_1, \dots, f_n$ is a strong Chebyshev system and $\mathrm{card}(\mathrm{supp}\mu)<n+1$, then the determinants $\det A^j$ (see below) can not be zero, which is impossible.
Let $\mu_k:=\frac{\mu|_{S_k}}{\mu(S_k)}$ and
$$ p_{i,k}:=\int_{S_k}f_id\mu_k-c_i, \ws\ws\ws q_k:=\int_{S_k}U^{\mu}d\mu_k-w, \ws\ws k=1, \dots, n+1, \ws i=1, \dots, n.$$
Notice that $\sum_{k=1}^{n+1}\mu(S_k)\int_{S_k}f_id\mu_k=c_i$ and $\sum_{k=1}^{n+1}\mu(S_k)\int_{S_k}U^{\mu}d\mu_k=w$. That is the homogeneous linear equation system $Ab=0$ has a strongly positive solution, $b>>0$, that is $b_j>0$, $j=1, \dots, n+1$, and we can assume that $\sum_{j=1}^{n+1}b_j=1$, where $A=\left[\begin{array}{ccc} p_{1,1}& \dots & p_{1,n+1} \\ \vdots &\dots & \vdots \\p_{n,1}& \dots & p_{n,n+1} \\ q_{1}& \dots & q_{n+1} \end{array}\right]$. (That is $\det A=0$.) It ensures that $A_1b=e$, where $A_1=\left[\begin{array}{ccc} p_{1,1}& \dots & p_{1,n+1} \\ \vdots &\dots & \vdots \\ p_{n,1}& \dots & p_{n,n+1} \\ 1& \dots & 1 \end{array}\right]$ and $e=\left[\begin{array}{c}0\\ \vdots \\0 \\1\end{array}\right]$. We can assume that the sets $S_k$ are chosen such that $\det A_1 \neq 0$. Let us define nonconstant positive continuous functions $h_k$ on $S_k$ such that $0<h_k(x)<2$ and $\int_{S_k}h_kd\mu_k=1$, $k=1, \dots, n+1$. Let $a_{i,k}:=\int_{S_k}f_ih_kd\mu_k-c_i$; $a_{n+1,k}:=\int_{S_k}U^{\mu}h_kd\mu_k-w$, $i=1, \dots, n$, $k=1, \dots, n+1$. Let $\Theta:=(\Theta_1, \dots, \Theta_{n+1})$, such that $\Theta_i \in (0,1)$, $i=1, \dots, n+1$ and
$$A(\Theta):=$$ $$\left[\begin{array}{ccc}(2\Theta_1-1)a_{1,1}+2(1-\Theta_1) p_{1,1}& \dots & (2\Theta_{n+1}-1)a_{1,n+1}+2(1-\Theta_{n+1}) p_{1,n+1} \\ \vdots &\dots & \vdots \\(2\Theta_1-1)a_{n,1}+2(1-\Theta_1)  p_{n,1}& \dots & (2\Theta_{n+1}-1)a_{n,1}+2(1-\Theta_{n+1})  p_{n,n+1} \\(2\Theta_1-1)a_{n+1,1}+2(1-\Theta_1)  q_{1}& \dots & (2\Theta_{n+1}-1)a_{n+1,n+1}+2(1-\Theta_{n+1}) q_{n+1} \end{array}\right].$$
Since by the indirect assumption $U^{\mu}(x)$  cannot be equivalently equal to $w$ $\mu$-a.e., we can assume that $\det A(\Theta)\not\equiv 0$. By the previous computation for $\Theta_0=\left(\frac{1}{2}, \dots, \frac{1}{2}\right)$ the equation $A(\Theta_0)b(\Theta_0)=0$ has a strongly positive solution which is also the solution of the equation $A_1(\Theta_0)b(\Theta_0)\\=e$. Since $b_j(\Theta_0)>0$ for all $j$, $\Theta_0$ has a (small) neighborhood, $V$ such that for any $\Theta \in V$ $b_j(\Theta)>0$ for all $j$, and $\det A_1(\Theta)$ has also the same sign as $\det A_1(\Theta_0)$.  Let us denote by $A^j(\Theta)\in \mathbb{R}^{n \times n}$ the matrix generated from $A(\Theta)$ such that we omit the last row and the $j$-th column. $b_j(\Theta)=\frac{(-1)^{n+1+j}\det A^j(\Theta)}{\det A_1(\Theta)}$, and so with any $\Theta \in V$ $\sum_{j=1}^{n+1}A(\Theta)_{n+1,j}b_j(\Theta)=\frac{\det A(\Theta)}{\det A_1(\Theta)}$. It is enough to show that there is a $\Theta \in V$ for which $\frac{\det A(\Theta)}{\det A_1(\Theta)}<0$, that is it leads to a contradiction since $\mu(\theta)=\sum_{j=1}^{n+1}b_j(\Theta)\mu_j \in \mathcal{M}_{\mathcal{F}}(X))$ and $E(\mu(\theta),\mu)<E(\mu)$. By the assumptions $\det A(\Theta_0)=0$ but $\det A(\Theta)\not\equiv 0$. Since the determinant consists of products which contains only one element from every row, in every products the variables $\Theta_i$ are at most on the first power. That is considering the determinant of $A(\Theta)$ as a function of the variables $\Theta_i$, all the elements of the main diagonal of the Hessian of $A(\Theta)$ are zeros, so the trace of the Hessian is zero, and so it has both positive and negative eigenvalues for all $\Theta \in V$, that is in $\Theta_0$ $A(\Theta)$ cannot has a minimum/maximum, so we can choose a suitable  $\Theta \in V$.

\medskip

As a corollary of the previous section, we are in position to state the main theorem of this section. In \cite{cw} after the cutting plane algorithm a discretization algorithm is introduced on compact intervals, which reduces the semi-infinite problem to finite ones. By Chebyshev-constant method we extend this idea to any compact metrizable spaces.

For simplicity let us denote by
$$\mathcal{R}(X):=\mathcal{R}_{\Phi,g}^e(X), \ws \ws \mbox{and}\ws\ws \mathcal{R}^{(n)}(X):=\mathcal{R}_{\Phi,g,Z_n}^e(X),$$
and let us recall that
$$\mathcal{R}^{(n)}_m(X):=\{\mu\in \mathcal{R}^{(n)}(X): \mathrm{Card}(\mathrm{supp}\mu)\leq m\}.$$

\medskip

\begin{theorem}
Besides the assumptions of $(CQP_e)$ and Lemma 10, let us suppose that the kernel $k(x,y)$ fulfils the Frostman's maximum principle. Then
\begin{equation}(CQP_e) = \inf_{\mu\in \mathcal{R}(X)}\int_X\int_X k(x,y)d\mu(x)d\mu(y)=w(\mathcal{R}(X))\end{equation}
\begin{equation}=\lim_{n\to \infty}w(\mathcal{R}^{(n)}(X))\end{equation}
\begin{equation}=\lim_{n\to\infty}q(\mathcal{R}^{(n)}(X))=\lim_{n\to\infty}\inf_{\mu \in \mathcal{R}^{(n)}(X)}\sup_{\nu\in\mathrm{Ex}\mathcal{R}^{(n)}(X)}\int_X\int_X k(x,y)d\mu(x)d\nu(y)\end{equation}
\begin{equation}=\lim_{n\to\infty}\underline{q}(\mathcal{R}^{(n)}(X))=\lim_{n\to\infty}\sup_{\mu \in \mathcal{R}^{(n)}(X)}\inf_{\nu\in\mathrm{Ex}\mathcal{R}^{(n)}(X)}\int_X\int_X k(x,y)d\mu(x)d\nu(y)\end{equation}
\begin{equation}=\lim_{n\to\infty}M(\mathcal{R}^{(n)}(X))\end{equation}
\begin{equation}=\lim_{n\to\infty}\lim_{m\to\infty}\sup_{\mu \in \mathcal{R}_m^{(n)}(X)}\inf_{\nu\in\mathrm{Ex}\mathcal{R}^{(n)}(X)}\int_X\int_X k(x,y)d\mu(x)d\nu(y)\end{equation}
\begin{equation}=\lim_{n\to\infty}\lim_{m\to\infty}\underline{q}(\mathcal{R}^{(n)}_m(X),\mathrm{Ex}\mathcal{R}^{(n)}(X)).\end{equation}
\end{theorem}

\proof (19) is the cutting plane algorithm (Lemma 8).  Let us recall that $\mathcal{R}^{(n)}(X)$ is $w^*$-compact and obviously is convex that is we can apply the results above. (20) follows from Lemma 10 and Corollary 2. Corollary 1 ensures (21). We get (22) from Theorem 3. (23) and (24) are Definition 5 and Definition 1.

\medskip

Examples and numerical corollaries of this theorem are discussed in the next section.

\section{Appendix}

To make more complete the discussion, in this section we cite some examples from the literature for kernel functions which satisfy the maximum principle and for sets of measures with given extremal sets.

Let us recall that $k$ satisfies the Frostman's maximum principle if for every measure $\mu \in \mathcal{M}(X)$ of compact support $\sup_{x\in X}U^{\mu}(x)=\sup_{x\in \mathrm{supp}\mu}U^{\mu}(x)$.

The Frostman's maximum principle for Newtonian kernel was proved by M. A. Maria \cite{m} and kernels of order $\alpha \leq 2$ by O. Frostman \cite{f}, for Riesz kernels see \cite{la}, for logarithmic kernel cf. e.g. \cite[Theorem 3.3.4]{ra} and \cite[Theorem B]{ki}. For more general kernels by L. Carleson \cite[page 14]{c}. For continuous kernels it is proved in \cite{fana} and \cite{y}, that the equivalence of the Chebyshev constant and the energy entails the maximum principle. Here we cite the result of Carleson, which is useful for applications. For further results see e.g. \cite{y} and \cite{ki}.

Let $H(t):\mathbb{R}\to \mathbb{R}_{+}$ be a non-negative continuous increasing, convex function. Let $x\in \mathbb{R}^d$ and $r=\|x\|$ the Euclidean norm of $x$ and $\Lambda(r)$ be a fundamental solution of Laplace's equation
$$\Lambda(r)=\left\{\begin{array}{ll}\log\frac{1}{r}, \ws d=2 \\ r^{2-d}, \ws d>2.\end{array}\right.$$
Let
$$k(x,y):=k(\|x-y\|),\ws \mbox{where} \ws
k(r)=H(\Lambda(r))\ws \mbox{if} \ws \int_0^{A}k(r)r^{d-1}dr<\infty.$$

\medskip

\begin{lemma}\cite[Theorem 1]{c}
If the kernel $k$ fulfils the assumptions above, it also fulfils the Frostman's maximum principle.\end{lemma}

\medskip

Now we turn to the discussion of extremal sets. Recalling the notation
$$\mathcal{R}_{\Phi,g,Z_n}^e(X)=\mathcal{M}_{\mathcal{F}}(X):=\{\mu \in \mathcal{M}(X) : \int_Xf_id\mu= c_i,\ws i=0,1,\dots, n\}.$$

R. Douglas gave descriptions of the extremal sets of sets of measures which fulfil some equality type conditions and which are contained by $\mathcal{M}(X)$. For more details see \cite{do} and \cite{do2}.

From the results of Douglas A. Karr derived a characterization of the extreme points of sets of measures of type $\mathcal{M}_{\mathcal{F}}(X)$.

\begin{lemma} \cite[Theorem 2.1, Proposition 2.1]{k}
Let $X$ be a compact, metrizable space.
For each $\mu\in \mathcal{M}_{\mathcal{F}}(X)$ the following assertions are equivalent:

\noindent {\rm (1)} $\mu$ is an extreme point of $\mathcal{M}_{\mathcal{F}}(X)$

\noindent {\rm (2/a)} $\mathrm{card}(\mathrm{supp} \mu)\leq n+1$ and

\noindent {\rm (2/b)} if $\mathrm{supp}=\{x_1,\dots, x_k\}$ then the vectors \\$(1, f_1(x_1), \dots, f_n(x_1)), \dots, (1, f_1(x_k), \dots, f_n(x_k))$ are linearly independent.

Furthermore $\mathcal{M}_{\mathcal{F}}(X)$ is nonempty if and only if $(c_1, \dots, c_n)$ lies in the closed convex hull of $\{(f_1(x), \dots, f_n(x)): x\in X\}$.\end{lemma}

\medskip

First let us suppose that $\mathcal{M}_{\mathcal{F}}(X)$ is nonempty. It is convenient to assume that $\{1, f_1,\dots, f_n\}$ is linearly independent cf. \cite{k}.
We can omit assumption (2/b) if $\{1, f_1,\dots, f_n\}$ is a Chebyshev system (for the definition see e.g. \cite{eb}). At this point we have the following form of e.g. (20) and (24).

\medskip

\begin{corollary}
If $X$ is a compact, metrizable Hausdorff space, $k$ is l.s.c. and fulfils the Frostman's maximum principle, $\mathcal{F}\subset L^1_{\mu}$ $\forall \mu \in \mathcal{M}_1(X)$ and is a Chebyshev system, $\mathcal{F}_n=\{f_0, f_1, \dots, f_n\}\subset \mathcal{F}$, then
$$CQP_e=w(\mathcal{R}_{\Phi,g}(X))$$ \begin{equation}=\lim_{n\to\infty}\inf_{\mu, \int_Xf_id\mu=c_i, \atop i=0,1,\dots, n}\sup_{y_1,\dots, y_l\in X,\ws l\leq n+1, \ws  \alpha_j>0\atop \sum_{j=1}^{l}\alpha_jf_i(y_j)=c_i, i=0,1,\dots, n} \sum_{j=1}^{l}\alpha_j\int_Xk(x,y_j)d\mu(x)\end{equation}
$$=\lim_{n\to\infty}\lim_{m\to\infty}\sup_{x_1,\dots, x_m\in X,\ws \beta_k>0\atop \sum_{k=1}^{m}\beta_kf_i(x_k)=c_i, i=0,1,\dots, n}\inf_{y_1,\dots, y_l\in X,\ws l\leq n+1, \ws  \alpha_j>0\atop \sum_{j=1}^{l}\alpha_jf_i(y_j)=c_i, i=0,1,\dots, n} \sum_{j=1}^{l}\sum_{k=1}^{m}\alpha_j\beta_kk(x_k,y_j).$$ \end{corollary}
That is we can handle the infinite (or continuous) problem as the limits of semi-infinite and finite minimax problems.

\medskip

\noindent{\bf Examples.}
There are several examples on Chebyshev systems (cf. e.g. \cite{eb}). We specify some of them with initial function $f_0\equiv 1$. Let $X =[a,b]\subset (0,\infty)$ and $0=\lambda_0<\lambda_1<\dots <\lambda_n$, then $\left\{x^{\lambda_0}, x^{\lambda_1},\dots ,x^{\lambda_n}\right\}$, $\left\{\cosh\lambda_0t, \cosh\lambda_1t,\dots ,\cosh\lambda_nt\right\}$ are Chebyshev systems. If $[a,b]\subset (-\infty,\infty)$, with the same exponents $\left\{e^{\lambda_0t},\dots ,
e^{\lambda_nt}\right\}$ is also a Chebyshev system.

It is also clear that $\mathcal{F}=\{x^{i} : i=0,1,\dots, n\}$ is also a Chebyshev system on $[a,b]=[0,1]$. Denoting the $r$th order difference operator by $\Delta^r$ the following corollary can be derived.

\medskip

\begin{lemma} \cite[Proposition 3.2]{k}

\noindent \rm{(a)} $\mathcal{M}_{\mathcal{F}}(X)$ is nonempty if and only if $(-1)^r\Delta^rc_k\geq 0$ for each $r,k\leq N$.

\noindent \rm{(b)} $\mu\in \mathcal{M}_{\mathcal{F}}(X)$ is an extreme point if and only if $\mathrm{card}(\mathrm{supp} \mu)\leq N+1$.\end{lemma}

\medskip

Here we mention another example (for further examples see \cite{k}). Let $f_i(x)=e^{-\lambda_ix}$ on $[0,1]$, where  $0=\lambda_0<\lambda_1<\dots <\lambda_N$. Let $1=c_0\geq c_1\geq \dots \geq c_N>0$. Then

\medskip

\begin{lemma} \cite[Proposition 3.3]{k}

\noindent \rm{(a)} $\mathcal{M}_{\mathcal{F}}(X)$ is nonempty if and only if $(c_1, \dots, c_N)$ is in the convex hull of the curve $\Gamma(u):=(u, u^{\frac{\lambda_2}{\lambda_1}}, \dots, u^{\frac{\lambda_N}{\lambda_1}})$, where $u\in \left[e^{-\lambda_1},1\right]$.

\noindent \rm{(b)} $\mu\in \mathcal{M}_{\mathcal{F}}(X)$ is an extreme point if and only if $\mathrm{card}(\mathrm{supp} \mu)\leq N+1$.\end{lemma}

\medskip

Finally we show an application of Corollary 1, which can be interesting in itself.

\medskip

\noindent {\bf Example.}
Let $H=K=[0,1]$ (as subsets of a suitable compact subset of $\mathbb{C}$, say) For instance let  $k(x,y)=\log\frac{1}{|x-y|}$. Let $\{c_i\}_{i=1}^N$ such that $(-1)^r\Delta^rc_k\geq 0$ for each $r,k\leq N$, $0=\lambda_0<\lambda_1<\dots <\lambda_M$, $(d_1,\dots, d_M)\in \mathrm{co}\Gamma u$ as above. Then
$$\sup_{\mu \atop \int_0^1x^{i}d\mu=c_i, \ws i=0,\dots, N}\inf_{\sum_{k=0}^l\alpha_k\delta_{y_k}, \ws y_k \in [0,1] \alpha_k>0, \ws l\leq M \atop \sum_{k=0}^l\alpha_ke^{-\lambda_jy_k}=d_j, \ws j=0,\dots, M}\sum_{k=0}^l\alpha_k\int_0^1 k(x,y_k)d\mu(x)$$ $$=\inf_{\nu \atop \int_0^1e^{-\lambda_jy}d\nu(y)=d_j,\ws j=0\dots, M}\sup_{\sum_{r=0}^s\beta_r\delta_{x_r} \ws x_r \in [0,1] \beta_r>0, \ws s\leq N\atop \sum_{r=0}^s\beta_rx_r^{i}=c_i, \ws i=0,\dots, N}\sum_{r=0}^s\beta_r\int_0^1 k(x_r,y)d\nu(y).$$

\medskip

\noindent \small{Department of Analysis, \newline
Budapest University of Technology and Economics}\newline
\small{ g.horvath.agota@renyi.mta.hu}


\medskip

\begin{thebibliography}{99}

\bibitem {ah} D. Adams, L. I. Hedberg, Function Spaces and Potential Theory, Springer-Verlag (1996)
\bibitem {a} E. J. Anderson, A Review of Duality Theory for Linear Programming over Topological Vector Spaces, {\it J. of Math. Analysis and Appl.} {\bf 97} 380-392, (1983)
\bibitem {eb} P. Borwein, T. Erd\'elyi, Polynomials and Polynomial Inequalities, Springer-Verlag (1996)
\bibitem {c} L. Carleson, Selected Problems on Exceptional Sets, D. Van Nostrand Co. Inc. (1967)
\bibitem{cw} S-Y. Chen, S-Y Wu, Algorithms for infinite quadratic programming in $L_p$ spaces, {\it J. of Comput. and Applied Math.} {\bf 213} 408-422, (2008)
\bibitem {ch} G. Choquet, Theory of Capacities, {\it Ann. Inst. Fourier} {\bf 5} Grenoble, 131-295 (1955)
\bibitem {do} R. G. Douglas, On extremal measures and subspace density, {\it Michigan Math. J.} {\bf 11} 243-246 (1964)
\bibitem {do2} R. G. Douglas, On extremal measures and subspace density II, {\it Proc. of the Amer. Math. Soc.} {\bf 17} (6) 1363-1365, (1966)
\bibitem {d} N. Dinculeanu, Integration on locally compact spaces, Noordhoff International Publishing Leyden (1974)
\bibitem {fana} B. Farkas, B. Nagy, Transfinite Diameter, Chebyshev Constant and Energy on Locally Compact Spaces, {\it Potential Anal} {\bf 28} 241-260, (2008)
\bibitem {fare} B. Farkas, Sz. R\'ev\'esz, Potential theoretic approach to rendezvous numbers, {\it Monatsh. Math.}, {\bf 148} (4), 309–331, (2006)
\bibitem {f} O. Frostman, Potentiel d'\'equilibre et capacit\'e des ensembles, {\it Comm. S\'em. Math. Univ. Lund} {\bf 3} (1935)
\bibitem {fu} B. Fuglede, On the Theory of Potentials in Locally Compact Spaces, {\it Acta Math.} {\bf 103}  140-215, (1960)
\bibitem {fu1} B. Fuglede, Le th\'eoreme du minimax et la th\'eorie fine du potentiel, {\it Ann. Inst. Fourier} {\bf 15} 65-87, (1965)
\bibitem {ho1} \'A. P. Horv\'ath, p-transfinite Diameter and p-Chebyshev Constant in Locally Compact Space, {\it Ann. Acad. Sci. Fenn. Math.} {\bf 40}, 851-874, (2015)
\bibitem {k} A. F. Karr, Extreme point of Certain Sets of Probability Measures, with Applications, {\it Math. of Operation Research} {\bf 8} (1) 74-85, (1983)
\bibitem {ki} M. Kishi, Unicity Principles in the Potential Theory, {\it Osaka Math. J.} {\bf 13}  41-74, (1961)
\bibitem {la} N. S. Landkof, Foundation of Modern Potential Theory, Springer-Verlag, New York, (1972)
\bibitem {m} M. A. Maria, The Potential of a Positive Mass and the Weight Function of Wiener, {\it Proc. Nat. Acad. Sci. USA} {\bf  20} 485-489, (1934)
\bibitem {ny} T. Nakamura, M. Yamasaki, Sufficient Conditions for Duality Theorems in Infinite Linear Programming Problems, {\it Hiroshima Math. J.} {\bf 9} 323-334 (1979)
\bibitem {o} M. Ohtsuka, On Potentials in Locally Compact Spaces, {\it J. Sci. Hiroshima Univ. ser. A-1.}  135-352, (1961)
\bibitem {o1} M. Ohtsuka, On Various Definitions of Capacity and Related Notions, {\it Nagoya Math.J.} {\bf 30} 121-127, (1967)
\bibitem {o2} M. Ohtsuka, A Generalization of Duality Theorem in the Theory of Linear Programming, {\it J. Sci. Hiroshima Univ. Ser. A-I} {\bf 30} 31-39, (1966)
\bibitem {o3} M. Ohtsuka, Generalized Capacity and Duality Theorem in Linear Programming, {\it J. Sci. Hiroshima Univ. Ser. A-I} {\bf 30} 45-56, (1966)
\bibitem {rrr} Q. Rajon, T. Ransford, J. Rostand, Computation of capacity via quadratic programming, {\it J. Math. Pures Appl.} {\bf 94} 398-413, (2010)
\bibitem {ra} T. Ransford, Potential Theory in the Complex Plane,  Cambridge Univ. Press (1995)
\bibitem {rh} J.Rigoberto Gabriel, O. Hern\'andez-Lerma, Strong Duality of the General Capacity Problem in Metric Spaces, {\it Math. Meth. Oper. Res.} {\bf 53} 25-34, (2001)
\bibitem {st} E. B. Saff, V. Totik, Logarithmic Potentials with External Fields, Springer-Verlag, (1997)
\bibitem {sch} H.H. Schaefer, M.P. Wolff, Topological Vector Spaces, Springer-Verlag New York, (1999)
\bibitem {y} M. Yamasaki, Relations between Capacities and Maximum Principle, {\it J. Sci. Hiroshima Univ. ser. A-1}  59-69, (1969)
\bibitem {yo} M. Yoshida, Some Examples Related to Duality Theorem in Linear Programming, {\it J. Sci. Hiroshima Univ. Ser. A-I} {\bf 30} 41-43, (1966)
\bibitem {w} S. Y. Wu, A cutting plane approach to solving quadratic infinite programs on measure spaces, {\it J. of Global Optimization} {\bf 21} 67-87, (2001)
\bibitem {wfl} S. Y. Wu, S. C. Fang, C. J. Lin, Solving General Capacity Problem by Relaxed Cutting Plane Approach, {\it Ann. of Oper. Res.} {\bf 103} 193–211, (2001)
\bibitem {wwt} Z. Wan, S. Y. Wu, K. L. Teo, Some properties on quadratic infinite programs of integral type, {\it Appl. Math. Letters} {\bf 20} 676–680, (2007)
\end{thebibliography}
\end{document}